\documentclass[final,1p,times]{elsarticle}

\usepackage{enumitem}
  \setlist{nosep} %

\usepackage[makeroom]{cancel}
\usepackage{amssymb,amsthm,mathtools,xfrac}
\newtheorem{remark}{Remark}[section]
\usepackage{stmaryrd}
\usepackage{verbatim}
\usepackage{soul}
\usepackage{comment}
\usepackage{color}
\usepackage{cleveref}
\usepackage{subcaption}
\renewcommand{\d}{\textnormal{d}}

\usepackage{siunitx}
\DeclareSIUnit\str{str} %
\DeclareSIUnit\sh{sh} %
\DeclareSIUnit\erg{erg} %

\usepackage[mathlines]{lineno}
\newcommand*\patchAmsMathEnvironmentForLineno[1]{%
  \expandafter\let\csname old#1\expandafter\endcsname\csname #1\endcsname
  \expandafter\let\csname oldend#1\expandafter\endcsname\csname end#1\endcsname
  \renewenvironment{#1}%
     {\linenomath\csname old#1\endcsname}%
     {\csname oldend#1\endcsname\endlinenomath}}%
\newcommand*\patchBothAmsMathEnvironmentsForLineno[1]{%
  \patchAmsMathEnvironmentForLineno{#1}%
  \patchAmsMathEnvironmentForLineno{#1*}}%
\AtBeginDocument{%
\patchBothAmsMathEnvironmentsForLineno{equation}%
\patchBothAmsMathEnvironmentsForLineno{align}%
\patchBothAmsMathEnvironmentsForLineno{flalign}%
\patchBothAmsMathEnvironmentsForLineno{alignat}%
\patchBothAmsMathEnvironmentsForLineno{gather}%
\patchBothAmsMathEnvironmentsForLineno{multline}%
}

\usepackage{xcolor}

\newcommand{\tcb}[1]{\textcolor{black}{#1}}

\usepackage{algorithm,algorithmicx}
\usepackage[noend]{algpseudocode}
\newcommand*\Let[2]{\State #1 $\gets$ #2}
\newcommand{\pluseq}{\mathrel{+}=}
\newcommand{\bOmega}{\boldsymbol{\Omega}}

\begin{document}

\begin{frontmatter}

\title{One-sweep moment-based semi-implicit-explicit\\integration for gray thermal radiation transport}
\author[LANL]{Ben S. Southworth}
\cortext[mycorrespondingauthor]{Corresponding author}
\ead{southworth@lanl.gov}
\author[LANL2]{Samuel Olivier}
\author[LANL]{H.K. Park}
\author[Tulane]{Tommaso Buvoli}

\address[LANL]{Theoretical Division, Los Alamos National Laboratory, P.O. Box 1663, Los Alamos, NM 87545 U.S.}
\address[LANL2]{Computer, Computational, and Statistical Sciences Division, Los Alamos National Laboratory, P.O. Box 1663, Los Alamos, NM 87545 U.S.}
\address[Tulane]{Department of Mathematics, Tulane University, 6823 St. Charles Avenue, New Orleans, LA 70118 US}

\begin{abstract}
Thermal radiation transport (TRT) is a time dependent, high dimensional partial integro-differential equation. In practical applications such as inertial confinement fusion, TRT is coupled to other physics such as hydrodynamics, plasmas, etc., and the timescales one is interested in capturing are often much slower than the radiation timescale. As a result, TRT is treated implicitly, and due to its stiffness and high dimensionality, is often a dominant computational cost in multiphysics simulations. Here we develop a new approach for implicit-explicit (IMEX) integration of gray TRT in the deterministic S$_N$ setting, which requires only one sweep per stage, with the simplest first-order method requiring only one sweep per time step. The partitioning of equations is done via a moment-based high-order low-order formulation of TRT, where the streaming operator and first two moments are used to capture the asymptotic stiff regimes of the streaming limit and diffusion limit. Absorption-reemission is treated explicitly, and although stiff, is sufficiently damped by the implicit solve that we achieve stable accurate time integration without incorporating the coupling of the high order and low order equations implicitly. Due to nonlinear coupling of the high-order and low-order equations through temperature-dependent opacities, to facilitate IMEX partitioning and higher-order methods, we use a semi-implicit integration approach amenable to nonlinear partitions. Results are demonstrated on thick Marshak and crooked pipe benchmark problems, demonstrating orders of magnitude improvement in accuracy and wallclock compared with the standard first-order implicit integration typically used.
\end{abstract}

\end{frontmatter}

\section{Introduction}

Consider the time-dependent grey thermal radiation transport (TRT) equations for angular intensity $I$ and temperature $T$ of the form:
\begin{subequations}\label{eq:ho}
\begin{align}\label{eq:ho-I}
  \frac{1}{c} \frac{\partial I}{\partial t} & = - \bOmega \cdot \nabla I - \sigma_t(T) I + \tfrac{1}{4\pi}ac\sigma_a(T)T^4,\\
  \rho c_v \frac{\partial T}{\partial t} & = - \sigma_a(T) acT^4 + \sigma_a(T)\int I\d\Omega ,\label{eq:ho-T}
\end{align}
\end{subequations}
where $c$ is the speed of light, $c_v$ is the specific heat (which in certain cases may depend on $T$), $\rho$ is the material density, and $\sigma_t(T)$ and $\sigma_a(T)$ are total and absorption opacities, respectively. In general, opacities are nonlinear functions of temperature and density, and often tabulated. Note that $I$ is high-dimensional, depending on time, space, and angle due to the collisional integral over direction of transport, $\int I\d\Omega$, and the coupled set of equations is very stiff due to the speed of light scaling. In practice, we often want to step well over the transport timescale (i.e., advection, $c\Delta t/\Delta x$, and absorption-emission, $\rho c_v/4ac\sigma T^3$ ), and as a result, in large-scale simulations such as coupled non-relativistic radiation hydrodynamics, \eqref{eq:ho} must be treated implicitly in time. Due to the high dimensionality and stiff behavior, implicit integration of \eqref{eq:ho} is very expensive, and often the bottleneck in multiphysics simulations. 

In this work, we use the discrete ordinate ($S_N$) method \cite{Adams.2002} to discretize the angular variable. A naive way of solving (\ref{eq:ho}) is via a Picard iteration of emission source,
\begin{subequations}\label{eq:picard}
\begin{align} \label{eq:picard-I}
    \frac{I^{n+1}_k - I^n}{c\Delta t_n} &+ \bOmega \cdot \nabla I^{n+1}_k + \sigma_t I^{n+1}_k = \frac{1}{4\pi} \sigma_a ac (T^{n+1}_{k-1})^4,\\
    \rho c_v\frac{T^{n+1}_k-T^n}{\Delta t_n} &= \sigma_a \int I^{n+1}_k d\Omega -\sigma_a ac (T^{n+1}_{k})^4 ,\label{eq:picard-T}
\end{align}
\end{subequations}
\tcb{where superscripts denote the time-step, subscripts denote the $k$th nonlinear Picard iteration, and $\Delta t_n = t_{n+1}-t_n$.}
Combination of (i) $S_N$ angular discretization, (ii) an upwind spatial discretization, and (iii) a Picard iteration of the emission source results in a block lower triangular matrix for $I$. Each angle couples through absorption-emission physics described in Eq.\ (\ref{eq:ho-T}). The process of inverting the block lower triangular matrix is often called a ``transport-sweep,'' and it is a key component to any efficient implicit solution strategy for S$_N$ transport. Due to the high dimensionality of $I$ and stiffness of absorption-emission coupling, a Picard iteration (i.e., source iteration) (\ref{eq:picard}) is computationally very expensive even with highly efficient transport sweep algorithms. To accelerate the convergence of stiff absorption-emission physics, a transport sweep is typically coupled with either a diffusion approximation \cite{Alcouffe.1977,Kopp.1963, Olson.2000} in temperature or a set of reduced moment equations \cite{gol1964quasi, Anistratov.1993}. The implicit iteration then alternates between a transport sweep and diffusion or moment-equation solve, repeating until convergence. 

In this paper, we propose a semi-implicit-explicit time integration scheme for TRT that requires only a single transport sweep per time step (or stage for higher order methods), and is stable while stepping over transport and collisional timescales. This is accomplished via a moment-based \cite{gol1964quasi,Anistratov.1993} high-order low-order (HOLO) \cite{Park.2012,Park.2013,Park.2020} formulation, where \eqref{eq:ho} is augmented with discretely-consistent LO moment equations to capture certain stiff behavior. We then decouple the HO system from the LO system by treating the (HO) emission source and opacities explicitly. This implicit-explicit (IMEX) linearization enables our algorithm to require a single HO transport sweep and implicit solution of nonlinear LO system per stage. We demonstrate that the proposed formulation can be significantly cheaper than fully implicit methods, while providing comparable accuracy and stability. Moreover, the proposed framework also allows for higher-order accuracy analogous to \cite{maginot2016high,LOU2019258}, without implicitly treating opacity, a significant challenge and cost, particularly when opacities are defined via tabular lookup or coupled to other physics such as hydrodynamics. \tcb{Although in this paper we focus specifically on the HOLO class of methods \cite{Park.2012,Park.2013,Park.2020}, the proposed IMEX framework naturally applies to other moment-based methods as well.}

IMEX methods have been proposed for transport and kinetic-type equations previously. \cite{McClarren.2008} proposes a ``semi-implicit'' integration for P$_N$ transport, which treats the collisional term implicitly, but the advection is still treated explicitly. This still poses significant stability constraints on time-step size in nonrelativistic settings such as radiation hydrodynamics, where one wants to closely follow the dynamical time scale, based on material velocity, which can be much larger than the radiation advection time scale based on the speed of light. \cite{Boscarino.2013} considers the stiff diffusion limit of hyperbolic equations plus relaxation, developing IMEX schemes that limit to an implicit method in the diffusion limit, but again remains explicit in the hyperbolic (transport) component. \cite{Dimarco.2013o9} develops asymptotic preserving (AP) IMEX methods for kinetic equations of Boltzmann type, limiting to an explicit scheme applied to the Euler equations as relaxation $\epsilon\to 0$. Similar asymptotic preserving analysis for the BGK equations is performed in \cite{Hu.2018}. In \cite{Chu.2019}, an IMEX scheme for relativistic neutrino transport in an astrophysics context is developed. There, the explicit advection time scale for radiation is not problematic, as that is the dynamical timescale of interest, but the local collision physics is stiffer and is treated using a local implicit nonlinear solve. A similar local implicit solve in the context of relativistic radiation hydrodynamics is considered in \cite{Weih.2020}, and related to the IMEX moment approximation coupled to hydrodynamics in \cite{Just.2015}. Although each of these works have proved useful in their respective area, none address the regime where $\Delta t \gg \Delta x/c$, and we cannot treat the transport terms explicitly.

The paper proceeds as follows. \Cref{sec:holo} introduces the HOLO formulation and finite-volume discretization for gray TRT, including important nuances on how we evaluate certain terms for higher order integration. Implicit-explicit and semi-implicit-explicit time integration is then introduced in \Cref{sec:time}. Our nonlinear partitioning of gray TRT is then presented in \Cref{sec:imex}, including details on how the algorithm is performed efficiently in practice. Numerical results demonstrating the efficacy of our new approach on stiff Marshak and crooked pipe TRT benchmark problems are provided in \Cref{sec:results}. In terms of both accuracy and efficiency, the new method is able to achieve orders of magnitude improvement in accuracy and/or wallclock time when compared with first-order implicit integration commonly used in practice. Conclusions and future work can be found in \Cref{sec:conc}.

\section{High-order low-order formulation}\label{sec:holo}

\subsection{Moment-based methods}

The HOLO algorithm is an efficient and robust iterative scheme for solving kinetic equations. The HOLO algorithm accelerates the solution of kinetic equations by capturing stiff collisional physics via a discretely-consistent, low-dimensional continuum description of the transport equation \cite{Park.2012,Park.2013,Park.2020}. 
Discrete closures are applied to the low-order system to compensate for physics and discretization discrepancies between the kinetic (HO) and continuum (LO) systems. 
For TRT, computational efficiency is achieved by implicitly treating the absorption-emission physics with the low-order system. 
This serves to isolate the high-dimensional transport equation from nonlinear material energy balance iterations. 

Following \cite{Park.2020}, we use lumped, linear discontinuous Galerkin discretizations of the S$_N$ transport equation as the high-order system. The LO system is constructed based on the following $P_1$ equations:
\begin{subequations}\label{eq:moment-eqs1}
\begin{align}
    \frac{\partial E}{\partial t} &= -\nabla \cdot \mathbf{F} - c\sigma_a(T)E + ac\sigma_a(T) T^4,\\
    \frac{1}{c}\frac{\partial \mathbf{F}}{\partial t} & = -\frac{c}{3}\nabla E - \sigma_t(T) \mathbf{F}, \label{eq:moment-eqs-F} 
\end{align}
\end{subequations}
where 
\begin{equation}\label{eq:moments}
  E \coloneqq \frac{1}{c}\int I\d\Omega,\\
  \hspace{5ex}
  \mathbf{F}\coloneqq \int \bOmega I\d\Omega,
\end{equation}
are the radiation energy density and radiative flux, respectively. 
The $P_1$ system is derived by approximating $\int \bOmega^2 I d\Omega \approx cE/3$. The HOLO algorithm modifies the spatially discretized $P_1$ equation in order to account for the physics fidelity of the high-order system and discretization mismatch by adding a discrete residual term into the radiative flux equation. The residual is evaluated by substituting the moments of HO solutions into the discretized first moment equations. 
This residual is normalized by the high-order radiation energy density. 
Let $\gamma(I)$ denote this normalized, discrete residual. Then HOLO iteratively solves \eqref{eq:ho-I} coupled to the following LO system: 
\begin{subequations}\label{eq:moment-eqs3}
\begin{align}
    \frac{\partial E}{\partial t} &= -\nabla \cdot \mathbf{F} - c\sigma_a(T)E + ac\sigma_a(T) T^4,\\
    \frac{1}{c}\frac{\partial \mathbf{F}}{\partial t} & = -\frac{c}{3}\nabla E - \sigma_t(T) \mathbf{F} + \gamma(I) cE, \label{eq:moment-eqs-F2} \\
    \rho c_v \frac{\partial T}{\partial t} & = -ac\sigma_a(T)T^4 + c\sigma_a(T) E.
\end{align}
\end{subequations}
The traditional HOLO algorithm applies an implicit backward Euler finite difference discretization in time (with explicitly linearized opacities). 
Given a Planck emission source, a transport sweep is performed to compute the intensity needed to form $\gamma$. 
The low-order system (\ref{eq:moment-eqs3}) is then iterated until convergence, yielding a new emission source for the transport equation. 
This algorithm requires a nested iteration where the outer iteration converges the coupling between the high and low-order systems and the inner iteration converges the nonlinear coupling between the low-order system and the material energy balance equation. 
Iterations are repeated until the high and low-order systems are consistent with each other. 

Note that this particular HOLO scheme is a member of a broader class of moment-based methods for kinetic equations (see \cite{CHACON201721} for a review) that includes the variable Eddington factor \cite{mihalas}, or quasidiffusion \cite{gol1964quasi}, method and the second moment method \cite{lewis_miller}. 
Such schemes differ from \cite{Park.2020} only in the design of the low-order system, and in principle the methods developed in this paper can be extended to other moment schemes. 

\subsection{$\gamma$ discretization}

In this paper, we define \emph{consistency} as the property that upon exact temporal integration of the discretized HO and LO systems, moments of the HO solution \eqref{eq:moments} match the LO solution. Note, this is slightly different from some existing literature which considers either steady problems or consistency following a single step of backward Euler \cite{Park.2012,Park.2013,Park.2020}. However, the underlying principle is analogous, and separating the consistency of HO and LO spatial discretizations from the temporal integration facilitates the development of consistent moment-based IMEX methods. Consistency of the moment systems ensures physics fidelity of the converged solution, and in our case we find is useful for the stability and convergence of the proposed IMEX integration scheme.

\begin{figure}[!htb]
\centering
    \includegraphics[width=0.75\textwidth]{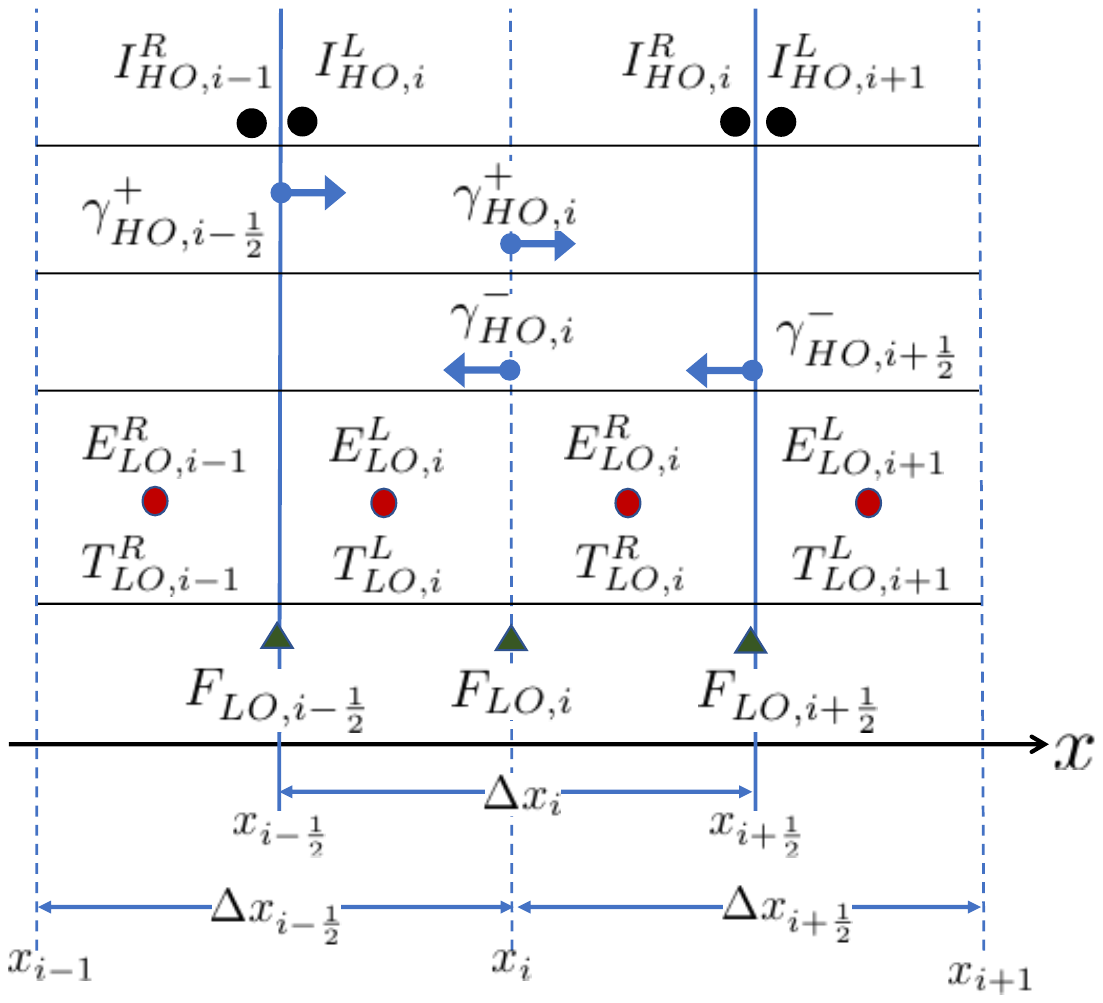}
\caption{Location of HO and LO variables.}
\label{fig:holo-varloc}
\end{figure}

In order to preserve the asymptotic diffusion limit, the LO system is discretized with (staggered) subcell finite volume method \cite{Park.2020}. For brevity, we present the discretization and equations in 1d, \tcb{with $\bOmega\cdot \nabla \mapsto \mu \tfrac{\partial}{\partial x}$}, but the framework and method extends to 2d or 3d without any technical complications. Define $\theta = aT^4$; then the discrete form of the moment equations from \eqref{eq:moment-eqs1} in 1d take the form
\begin{subequations}\label{eq:moment-eqs2}
\begin{align}
    \frac{\Delta x_i}{2}\frac{\partial E_{LO,i}^L}{\partial t} + 
        \left[{F}_{LO,i}-F_{LO,i-\frac{1}{2}}\right] + \frac{\sigma_i\Delta x_i}{2}cE_{LO,i}^L & = \frac{\sigma_ic\Delta x_i\theta_{LO,i}^L}{2}, \\
    \frac{\Delta x_i}{2}\frac{\partial E_{LO,i}^R}{\partial t} +
        \left[{F}_{LO,i+\frac{1}{2}}-F_{LO,i}\right] + \frac{\sigma_i\Delta x_i}{2}cE_{LO,i}^R & = \frac{\sigma_ic\Delta x_i\theta_{LO,i}^R}{2}, \\
    \frac{1}{c}\frac{\partial {F}_{LO,i-\frac{1}{2}}}{\partial t} +
        \frac{c}{3}\frac{E_{LO,i}^L - E_{LO,i-1}^R}{\Delta x_{i-\frac{1}{2}}/2} + \sigma_{i-\frac{1}{2}}{F}_{LO,i-\frac{1}{2}} & = \gamma^+_{HO,i-\frac{1}{2}}cE_{LO,i-1}^R-\gamma^-_{HO,i-\frac{1}{2}}cE_{LO,i}^L, \label{eq:F-face}\\
    \frac{1}{c}\frac{\partial {F}_{LO,i}}{\partial t} +
        \frac{c}{3}\frac{E_{LO,i}^R - E_{LO,i}^L}{\Delta x_{i}/2} + \sigma_{i}{F}_{LO,i} & = \gamma^+_{HO,i}cE_{LO,i}^L-\gamma^-_{HO,i}cE_{LO,i}^R, \label{eq:F-interior}
\end{align}
\end{subequations}
where superscripts $L,~R$, and $\pm$ denote left and right states in cell $i$ and positive ($\mu>0$) and negative ($\mu<0$) angular upwinding directions, respectively. The location of variables are shown in Fig.\ \ref{fig:holo-varloc}.

In order to define equations for $\gamma$, we use the partial radiative fluxes:
\begin{align}\label{eq:partial-F}
    \begin{split}
        {F}_{HO,i}^+ & = \sum_{\mu_m>0} \mu_m\omega_m\frac{I_{m,i}^L + I_{m,i}^R}{2}, \\
        {F}_{HO,i}^- & = \sum_{\mu_m<0} |\mu_m|\omega_m\frac{I_{m,i}^L + I_{m,i}^R}{2},
    \end{split}
    \hspace{-4ex}
    \begin{split}
        \hat{{F}}_{HO,i-\frac{1}{2}}^+ & = \sum_{\mu_m>0} \mu_m\omega_m\hat{I}_{m,i-\frac{1}{2}}, \\
        \hat{{F}}_{HO,i-\frac{1}{2}}^- & = \sum_{\mu_m<0} |\mu_m|\omega_m\hat{I}_{m,i-\frac{1}{2}}, 
    \end{split}
\end{align}
for face $i-\tfrac{1}{2}$ and cell interior $i$. $\omega$ is an angular quadrature weight. From here, we define our four consistency values $\{\gamma^+_{HO,i-\frac{1}{2}}$, $\gamma^+_{HO,i}$, $\gamma^-_{HO,i-\frac{1}{2}}$, $\gamma^-_{HO,i}\}$ to satisfy the following auxiliary equations:
\begin{subequations}\label{eq:aux-eqs}
\begin{align}
    cE_{HO,i}^L \gamma^+_{HO,i} & = 
        \frac{1}{c}\frac{\partial{{F}}^+_{HO,i}}{\partial t} + \frac{c}{3}\frac{E_{HO,i}^R - E_{HO,i}^L}{\Delta x_{i}} + \sigma_{i}{{F}}^+_{HO,i}, \label{eq:gamma1}\\
    cE_{HO,i}^R \gamma^-_{HO,i} & = 
        \frac{1}{c}\frac{\partial{{F}}^-_{HO,i}}{\partial t} - \frac{c}{3}\frac{E_{HO,i}^R - E_{HO,i}^L}{\Delta x_{i}} + \sigma_{i}{{F}}^-_{HO,i}, \label{eq:gamma2}\\
    cE_{HO,i-1}^R \gamma^+_{HO,i-\frac{1}{2}} & = 
        \frac{1}{c}\frac{\partial\hat{{F}}^+_{HO,i-\frac{1}{2}}}{\partial t} + \frac{c}{3}\frac{E_{HO,i}^L - E_{HO,i-1}^R}{\Delta x_{i-\frac{1}{2}}} + \sigma_{i-\frac{1}{2}}\hat{{F}}^+_{HO,i-\frac{1}{2}}, \label{eq:gamma3}\\
    cE_{HO,i}^L \gamma^-_{HO,i-\frac{1}{2}} & = 
        \frac{1}{c}\frac{\partial\hat{{F}}^-_{HO,i-\frac{1}{2}}}{\partial t} - \frac{c}{3}\frac{E_{HO,i}^L - E_{HO,i-1}^R}{\Delta x_{i-\frac{1}{2}}} + \sigma_{i-\frac{1}{2}}\hat{{F}}^-_{HO,i-\frac{1}{2}}.\label{eq:gamma4}
\end{align}
\end{subequations}
Note that summing Eqs.\ (\ref{eq:gamma1}) and (\ref{eq:gamma2}), and Eqs.\ (\ref{eq:gamma3}) and (\ref{eq:gamma4}) yields (\ref{eq:F-interior}) and (\ref{eq:F-face}), respectively. 

While not discussed in \cite{Park.2020}, boundary conditions are applied by closing the radiative flux at the boundary of the domain. 
We consider two types of boundary conditions derived by manipulating partial fluxes at the boundary referred to as full and half range. 
First, we propose the LO boundary radiative flux $F_n^{B} = \mathbf{F}^{B}\cdot \mathbf{n}_{B}$ to take a form of:
\begin{align}
   F_{n,LO}^{B} = \gamma^{B}_{1} cE_{LO}- \gamma^{B}_{0},
\end{align}
where $\mathbf{n}_{B}$ is the outward normal vector to a boundary face, and $\gamma^{B}_{0,1}$ are the parameters which depends on the HO solutions. The boundary HO radiative flux is
\begin{align}
    F_{n,HO}^{B} &= \int_{\bOmega \cdot \mathbf{n}_B>0}\bOmega \cdot \mathbf{n}_{B} I d\Omega-\int_{\bOmega \cdot \mathbf{n}_B<0}|\bOmega \cdot \mathbf{n}_{B}| I d\Omega,\nonumber\\
    &= F_{HO}^{B,out} - F_{HO}^{B,in}.
\end{align}

Half range boundary conditions simply split the radiative flux into incoming and outgoing and directly apply HO incoming radiative flux, thus $\gamma_{0,1}^{B}$ becomes: 
\begin{align}
    \gamma_0^{B} &= F_{HO}^{B,in},\hspace{5ex}
    \gamma_1^{B} =\frac{ F^{B,out}_{HO}}{cE_{HO}}.
\end{align}
The full range boundary conditions manipulate the partial fluxes such that 
\begin{align}
      F =  F^{B,out}_{HO}-F^{B,in}_{HO}  = (F^{B,out}_{HO} + F^{B,in}_{HO})-2F^{B,in}_{HO}.
\end{align}
Then $\gamma_{0,1}^{B}$ becomes:
\begin{align}
    \gamma_0^{B} &= 2F_{HO}^{B,in},\hspace{5ex}
    \gamma_1^{B} =\frac{ F^{B,out}_{HO}+F_{HO}^{B,in}}{cE_{HO}} .
\end{align}
For a fully implicit scheme, the choice of boundary conditions along with the definition of the cell-edge opacity, $\sigma_{i-1/2}$, are inconsequential; the discrete residual corrects the low-order system such that it matches the high-order solution to the iterative tolerance regardless of the choice of boundary condition or interface opacity. 
However, with IMEX time integration, the high and low-order systems are no longer iterated to convergence introducing a temporal inconsistency between high and low-order. 
As a result, the boundary conditions and interface opacity treatments do effect the final solution, and were seen to significantly impact solution quality and accuracy. 
In \Cref{sec:results:marshak}, we compare the use of the maximum, minimum, and harmonic average of the opacities in neighboring cells for $\sigma_{i-1/2}$ and the use of half and full range boundary conditions. \tcb{Nevertheless, we note that these schemes are conservative in the moments regardless of convergence of the consistency terms.} For more details on HOLO and the discretization see \cite{Park.2020}.\footnote{Note, there is an index typo in \cite{Park.2020} that has been corrected in \eqref{eq:aux-eqs}}

\subsection{Evaluation of HOLO $\gamma$ consistency term}\label{sec:holo:gamma}

Note that the $\gamma$ coefficients in \eqref{eq:aux-eqs} are defined in semi-discrete form (i.e., continuous time derivative). When the HO equation is discretized with a simple finite-difference scheme such as backward Euler, the time derivative in \eqref{eq:aux-eqs}, e.g., $\partial\hat{\mathbf{F}}/\partial t$, can also be evaluated with the same finite-difference scheme \cite{Park.2020}. 

In this subsection, we derive an expression for a time derivative of HO moments in terms of RK-stage solutions of the HO system. This enables one to use \eqref{eq:aux-eqs} for general IMEX-RK methods, where we cannot have a $\partial/\partial t$ directly in the evaluation of $\gamma$. The time-derivative in \eqref{eq:aux-eqs} must be derived consistently from the HO system, e.g.,
\begin{align}
    \frac{\partial F_{HO}^{\mp}}{\partial t} =  \frac{\partial \int_{\mu\lessgtr 0} |\mu| I d\mu}{\partial t} =\int_{\mu\lessgtr 0} |\mu| \frac{\partial I}{\partial t} d\mu.
\end{align}
The time derivative for $I$ at each RK stage $k$ can be derived by rearranging the formula for the $k$th RK stage and using (\ref{eq:ho-I}):
\begin{align}\label{eq:ho-solve}
\begin{split}
    I_{(k)} & = \hat{I}_{(k)} + \Delta t \left(c\bOmega \cdot \nabla I_{(k)} - c\sigma_t(T) I_{(k)} + \tfrac{1}{4\pi}ac^2\sigma_a(T)T^4\right) \\
    & = \hat{I}_{(k)} + \Delta t\frac{\partial I}{\partial t}(I_{(k)},T)
\end{split}
\end{align}
As we discuss in the Section \ref{sec:time}, each implicit RK stage solves (\ref{eq:ho-solve}) for $I_{(k)}$ with fixed $T$. $\hat{I}$ includes any forcing terms (in particular, artificial forcing arising from previous stage values). Further rearranging (\ref{eq:ho-solve}) yields the desired time-derivative expression for stage $k$: %
\begin{equation}\label{eq:dIdt}
    \frac{\partial I}{\partial t}(I_{(k)},T) = \frac{I_{(k)}-\hat{I}_{(k)}}{\Delta t}.
\end{equation}
Taking discrete moments of Eq.\ (\ref{eq:dIdt}) at interior and face points provides the information we need to compute $\gamma$. In particular, $E_{HO}$  in \eqref{eq:aux-eqs} is evaluated directly from the angular intensity solution $I_{(k)}$ following the transport sweep.
\begin{align}
    E_{HO,i}^{R,L} = \sum_m \omega_m \left[I_{(k)}\right]_{m,i}^{R,L}
\end{align}
On the other hand, to compute time-derivatives of partial radiative flux, we appeal to the definitions in \eqref{eq:partial-F} and \eqref{eq:dIdt}, e.g., at the $k$th time index,
\begin{subequations}\label{eq:indirect-partial-F}
\begin{align}
    \frac{\partial{F}_{HO,i}^+}{\partial t} & = \sum_{\mu_m>0} \frac{\mu_m\omega_m}{2} \left( \left[\tfrac{\partial I_{(k)}}{\partial t}\right]_{m,i}^L + \left[\tfrac{\partial I_{(k)}}{\partial t}\right]_{m,i}^R\right) \\
    & = \sum_{\mu_m>0} \frac{\mu_m\omega_m}{2} \left( \left[\tfrac{I_{(k)}-\hat{I}_{(k)}}{\Delta t}\right]_{m,i}^L + \left[\tfrac{I_{(k)}-\hat{I}_{(k)}}{\Delta t}\right]_{m,i}^R\right).
\end{align}
\end{subequations}
Similar computations yield $\partial{{F}}^-_{HO,i}/\partial t$, $\partial\hat{{F}}^+_{HO,i-\frac{1}{2}}/\partial t$, and $\partial\hat{{F}}^-_{HO,i-\frac{1}{2}}/\partial t$ to use in \eqref{eq:aux-eqs}, all as a function of the angular intensity variables $I_{(k)}$ and $\hat{I}_{(k)}$. 

\section{Semi-implicit-explicit time integration}\label{sec:time}

\subsection{Runge-Kutta methods}\label{sec:time:rk}

Consider the nonlinear autonomous set of ODEs
\begin{align}\label{eq:ode}
    \frac{\partial y}{\partial t} = \mathcal{N}(y).
\end{align}
Runge-Kutta (RK) methods \cite{butcher1996history} are a popular class of methods to approximate the integration of \eqref{eq:ode}. We adopt the standard notation of a Butcher tableaux,
$
\begin{array}
{c|c}
\mathbf{c} & A\\
\hline
&\mathbf{b}
\end{array},
$
where $A$ corresponds to coefficients, $\mathbf{b}$ to weights, and $\mathbf{c}$ to quadrature points within a time-step. The corresponding $s$-stage RK method takes the form
\begin{align*}
    Y_{(i)} &= y^n + \Delta t \sum_{j=1}^s a_{ij} \mathcal{N}(Y_{(j)})\hspace{4ex}\text{for }i=1,..,s,\\
    y^{n+1} &= y^n + \Delta t \sum_{j=1}^s b_j \mathcal{N}(Y_{(j)}).
\end{align*}
As mentioned previously, non relativistic transport codes almost universally consider some form of implicit integration. For implicit RK methods, we only consider diagonally-implicit RK (DIRK) methods here, which assume that $A$ is lower triangular. This makes for simpler implicit solves than fully implicit RK methods with dense $A$ and implicit coupling between all stages. 

Even restricting ourselves to DIRK methods, the implicit solution of the transport equations remains very computationally expensive due to the high dimensionality of the problem and stiff nonlinear coupling between variables. In this paper we seek an implicit-explicit-like method that is stable to take timesteps on the physical timescale that we are interested in, while also providing a cheaper evolution of the equations than purely implicit methods. Significant research has been done on partitioned and additive RK (ARK) methods \cite{Ascher.1997,Kennedy.2003tv4}, which provide a formal framework for integrating equations that have an additive partition into stiff and nonstiff parts, $\mathcal{N}(y) = \mathcal{N}_E(y) + \mathcal{N}_I(y)$, where we treat $\mathcal{N}_E$ explicitly and $\mathcal{N}_I$ implicitly. However, in realistic applications the equations do not always have a clear additive separation between implicit and explicit operators; e.g., see the recent paper on radiation hydrodynamics \cite{southworth2023implicit}. In our case, the transport equations have nonlinear coupling between variables that is difficult to separate in an additive fashion. The primary examples that arise in our partitioning of the transport equations (see \cref{sec:imex}) are opacities, e.g., $\sigma_t(T)I$ and the HOLO $\gamma$ correction. In each of these cases, to reduce the computational cost of implicit solves we do not want to treat the full multi-variable terms implicitly, but we also cannot treat these terms fully explicitly for stability reasons. 

\subsection{Semi-implicit RK methods}\label{sec:time:sirk}

To mitigate the nonlinear coupling issue, we follow the semi-implicit strategy of \cite{Boscarino.2015,Boscarino.2016}, where we introduce an auxiliary $*$-variable and consider IMEX-RK methods to solve the following coupled equations:
\begin{equation}\label{eq:mod-ode}
    \frac{\partial y^*}{\partial t} = \mathcal{N}(y^*,y), \hspace{3ex}
    \frac{\partial y}{\partial t} = \mathcal{N}(y^*,y).
\end{equation}
{\color{black} Note that because the right-hand sides of $y$ and $y^*$ are identical, given by $\mathcal{N}(y^*,y)$, integrating the coupled system exactly in time will yield identical solutions, $y(t) = y^*(t)$, which are also identical to exact integration of the original equation \eqref{eq:ode}. By duplicating the equations, however, we are able to apply discrete partitioned integrators to the coupled system \eqref{eq:ode-part}, which are explicit in $y^*$ and implicit in $y$, regardless of nonlinear interaction between the two variables. In doing this, we provide a flexibility and implementation very similar to standard Lie-Trotter operator splitting \emph{inside of the RK stage solutions}, thereby incurring the benefits of RK integration, and avoiding limitations of standard operator splitting. We will carefully formulate our problem so that we evaluate stiff terms implicitly at $y$ and non-stiff terms explicitly at $y^*$.} In this form, we can apply general IMEX-RK schemes to the modified formulation in \eqref{eq:mod-ode} while allowing nonlinear partitioning, such as $\sigma_t(T^*)I$, which will be treated explicitly in $T$ and implicitly in $I$. Define an IMEX-RK Butcher tableaux
\begin{equation}\label{eq:imex-tableaux}
\renewcommand\arraystretch{1.2}
\begin{array}
{c|c}
\tilde{\mathbf{c}} & \tilde{A}\\
\hline
& \mathbf{b}
\end{array}
,\hspace{3ex}
\begin{array}
{c|c}
\mathbf{c} & A\\
\hline
&\mathbf{b}
\end{array}
\end{equation}
where tilde-coefficients denote an explicit scheme ($\tilde{A}$ is strictly lower triangular, $\tilde{a}_{ij} = 0,~j \ge i$) and non-tilde coefficients represent a DIRK scheme ($A$ is lower triangular). Often, the time-quadrature points for the explicit and implicit schemes are assumed to be equal, $\tilde{\mathbf{c}} = \mathbf{c}$; however, since we are focused on autonomous problems, the quadrature points in time do not directly enter our integration scheme, so either case is equivalent from an implementation perspective. Note, we have assumed that $\mathbf{b} = \tilde{\mathbf{b}}$, \tcb{which leads to an updated solution $y_{n+1} = y_{n+1}^*$; thus even though there are distinct stages in the $*$- and non-$*$-variables, the actual solution updates are identical, eliminating the need to track two distinct discrete solutions or choose which one to keep.} With this choice, the resulting semi-implicit-explicit RK (SIMEX-RK) scheme takes the form
\begin{subequations}\label{eq:limex}
\begin{align}\label{eq:limex-exp}
    Y^*_{(i)} &= y_n + \Delta t \sum_{j=1}^{i-1} \tilde{a}_{ij} \mathcal{N}(Y^*_{(j)},Y_{(j)}),\\
    Y_{(i)}   &= y_n + \Delta t \sum_{j=1}^i {a_{ij}} \mathcal{N}(Y^*_{(j)},Y_{(j)}),\label{eq:limex-imp} \\
    y_{n+1} &= y_n + \Delta t \sum_{j=1}^s {b_{j}} \mathcal{N}(Y^*_{(j)},Y_{(j)}).
\end{align}
\end{subequations}
The simplest example of such a method is the first-order linearly implicit explicit Euler method, where
\begin{equation}\label{eq:limex-o1}
    y_{n+1} = y_n + \Delta t \mathcal{N}(y_n, y_{n+1}),
\end{equation}
which is effectively what has long been used in transport simulations to solve the equations implicitly with opacities evaluated at the previous time step \cite{Larsen.1988}.

In cases where there is a true additive partition of implicit and explicit operators, it can be more natural to separate the operator into an explicit and semi-implicit component:
\begin{align}\label{eq:ode-split2}
    \frac{\partial}{\partial t}
    y = \mathcal{N}_E(y^*) + \mathcal{N}_I(y^*,y)
\end{align}
Here $\mathcal{N}_E(y^*)$ is all explicit physics operators; $\mathcal{N}_I(y^*,y)$ is implicit, but potentially only in some components of $y$, and we leave this separation to the user.  A simple algorithm constructed that avoids evaluating the spatial discretization operators as much as possible is provided in \Cref{alg:limex}. Following \Cref{alg:limex}, the $j$th SIMEX-RK stage then takes the form
\begin{align}
    \mathbf{r}_{(j)} & = \mathbf{y}_n + \Delta t\sum_{k=1}^{j-1} \mathcal{N}_I(\mathbf{Y}^*_{(k)}, \mathbf{Y}_{(k)}), \\
    \mathbf{Y}_{(j)} - a_{jj}\Delta t\mathcal{N}_I(\mathbf{Y}^*_{(j)}, \mathbf{Y}_{(j)}) & = 
        \mathbf{r}_{(j)} + a_{jj}\Delta t\mathcal{N}_E(\mathbf{Y}^*_{(j)}), \\
    \mathcal{N}_I(\mathbf{Y}^*_{(j)}, \mathbf{Y}_{(j)}) & = \left(\mathbf{Y}_{(j)} -
        \mathbf{r}_{(j)} - a_{jj}\Delta t\mathcal{N}_E(\mathbf{Y}^*_{(j)})\right) / (a_{jj}\Delta t),
\end{align}
where we indirectly evaluate the implicit operator following the implicit solve. Note, here we expect the user-provided implicit solve to not only solve for the solution of the implicit variables, but also update any explicit variables/equations that depend on the implicit solution, such as absorption-reemission in rad-hydro. 
\begin{algorithm}[!htb]
  \caption{SIMEX-RK time step
    \label{alg:limex}}
  \begin{algorithmic}[1]
    \Let{$\mathbf{r}_{(i)}$}{$\mathbf{0}$ for $i=1,...,s$}
        \Comment{Initialize implicit right-hand side vectors}
    \Let{$\mathbf{Y}_{(i)}^*$}{$\mathbf{y}_n$ for $i=1,...,s$}
        \Comment{Initialize explicit stage vectors}
    \Let{$\mathbf{y}_{n+1}$}{$\mathbf{y}_n$}
        \Comment{Initialize solution}
    \Statex{}

    \For{$j = 1 \textrm{ to } s$}\Comment{Loop over stages $1,...,s$}

        \If{$a_{jj} \neq 0$}\Comment{Standard implicit stage}
            \State{$\boldsymbol{\delta} = \Delta t\mathcal{N}_E(\mathbf{Y}^*_{(j)})$}
            \Comment{Evaluate explicit part of operator}\label{line:delta}
            \State{$\mathbf{r}_{(j)} \pluseq \mathbf{y}_n + a_{jj}\boldsymbol{\delta}$}
            \Comment{Update implicit right-hand side}\label{line:rhs}
            \State{Solve $\mathbf{Y} - a_{jj}\Delta t \mathcal{N}_I(\mathbf{Y}^*_{(j)},\mathbf{Y}) = \mathbf{r}_{(j)}$}
            \Comment{Implicit solve for $\mathbf{Y}$}
            \State{$\boldsymbol{\delta} \pluseq (\mathbf{Y}-\mathbf{r}_{(j)})/a_{jj}$}\label{line:N}
            \Comment{$\boldsymbol{\delta} \mapsto \Delta t\mathcal{N}(\mathbf{Y}^*_{(j)},\mathbf{Y}_{(j)})$}
        \Else\Comment{Explicit stage in implicit variable}
            \State{$\mathbf{Y} = \mathbf{y}_n + \Delta t\mathbf{r}_{(j)}$} 
            \State{$\boldsymbol{\delta} = \Delta t\mathcal{N}(\mathbf{Y}^*_{(j)},\mathbf{Y})$}
        \EndIf
        \Statex
        \State{$\mathbf{y}_{n+1} \pluseq b_{j}\boldsymbol{\delta}$} 
            \Comment{Update solution with $j$th residual}
        \For{$i =j+1 \textrm{ to } s$}\label{line:loop}
            \State{$\mathbf{Y}_{(i)}^* \pluseq \tilde{a}_{i,j}\boldsymbol{\delta}$}
            \Comment{Update future explicit stages with $j$th residual}
            \State{$\mathbf{r}_{(i)} \pluseq {a}_{i,j}\boldsymbol{\delta}$}
            \Comment{Update future implicit right-hand sides with $j$th residual}
        \EndFor
    \EndFor
    \Statex
    \State\Return $\mathbf{y}_{n+1}$
    \end{algorithmic}
\end{algorithm}

{\color{black}
We have tested many IMEX-RK schemes for TRT; similar to our recent paper on IMEX-integration for gray radiation hydrodynamics \cite{southworth2023implicit}, we find the following schemes to be consistently superior in terms of robustness (stability) and accuracy:
\begin{itemize}
    \item \emph{H-LDIRK2(2,2,2)} \cite[Table II]{Pareschi.2005}, a two-stage 2nd-order method consisting of 2nd-order L-stable SDIRK implicit and 2nd-order SSP explicit methods:
    \begin{align*}
        \begin{array}{c | c c c}
        0 & 0 & 0 \\
        1 & 1 & 0 \\\hline
        & 1/2 & 1/2
        \end{array},
    \hspace{5ex}
        \begin{array}{c | c c}
        \gamma & \gamma & 0 \\
        1 -\gamma & 1-2\gamma & \gamma \\\hline
        & 1/2 & 1/2
        \end{array},
    \hspace{2ex} \textnormal{for }\gamma=1-1/\sqrt{2}.
    \end{align*}

	 \item \emph{SSP-LDIRK2(3,3,2)} \cite[Table IV]{Pareschi.2005}, a three-stage, 2nd-order method consisting of 2nd-order L-stable DIRK implicit and 2nd-order SSP explicit methods,
    \begin{align*}
        \begin{array}{c | c c c}
        0 & 0 & 0 & 0 \\
        1/2 & 1/2 & 0 & 0\\
        1 & 1/2 & 1/2 & 0 \\\hline
        & 1/3 & 1/3 & 1/3
        \end{array},
    \hspace{5ex}
        \begin{array}{c | c c c}
        1/4 & 1/4 & 0 & 0 \\
        1/4 & 0 & 1/4 & 0 \\
        1/3 & 1/3 & 1/3 & 1/3 \\\hline
        & 1/3 & 1/3 & 1/3
        \end{array}.
    \end{align*}

	 \item \emph{SSP-LDIRK3(3,3,2)} \cite[Table V]{Pareschi.2005}: a three-stage, 2nd-order method consisting of 2nd-order L-stable DIRK implicit and 3rd-order SSP explicit methods.
    \begin{align*}
        \begin{array}{c | c c c}
        0 & 0 & 0 & 0 \\
        1 & 1 & 0 & 0\\
        1/2 & 1/4 & 1/4 & 0 \\\hline
        & 1/6 & 1/6 & 2/3
        \end{array},
    \hspace{5ex}
        \begin{array}{c | c c c}
        \gamma & \gamma & 0 & 0 \\
        1-\gamma & 1-2\gamma & \gamma & 0 \\
        1/2 & 1/2 - \gamma & 0 & \gamma \\\hline
        & 1/6 & 1/6 & 2/3
        \end{array},
        \hspace{2ex} \textnormal{for }\gamma=1-1/\sqrt{2}.
    \end{align*}
\end{itemize} 
as well as the simplest case of 1-stage LIMEX-Euler \eqref{eq:limex-o1}. 
}

\begin{remark}[Nonlinearly partitioned RK methods]
We have begun developing a new class of RK methods specifically designed for nonlinearly partitioned problems, that is, the methods are applied directly to an equation of the form $y' = F(y,y)$ without any copying of variables/equations \cite{nprk1,nprk2}. \tcb{Here we are free to choose where we evaluate $F$ in the first or second argument, and we then develop NPRK methods that, for example, treat the first argument explicitly and the second argument implicitly. By developing these methods for arbitrary functions of two arguments, this naturally faciltiates nonlinearly partitioned integration.} We do not explore these methods in detail here, but plan to develop methods specifically for TRT in future work, and initial results demonstrate improved stability on stiff Marshak problems \cite{nprk1}.
\end{remark}

\section{Moment-based semi-implicit-explicit integration for TRT}\label{sec:imex}

Each implicit transport solve is typically very expensive due to the high dimensionality, ill-conditioning, and advective nature of the equations. In practice, it is common to not fully converge linear or nonlinear implicit iterations, instead doing one or a few iterations and moving on. Here we formalize such a strategy to allow for guaranteed accuracy and higher order methods, with one sweep per stage. To this end, we appeal to the auxiliary variable as in \eqref{eq:mod-ode} and write our complete (HO and LO) transport equation as\footnote{\tcb{Note, here \eqref{eq:limex-trt1} $I$ and $E$ indicate equations associated with state variables angular intensity and energy, rather than implicit and explicit partitions of a general operator as used previously \eqref{eq:ode-split2}.}}
\begin{subequations}\label{eq:limex-trt1}
\begin{align}
    \frac{\partial I}{\partial t}  &= \tfrac{1}{4\pi}ac^2\sigma_a(T^*) (T^*)^4 - c\bOmega \cdot \nabla I - c\sigma_t(T^*) I &&\coloneqq \mathbf{N}_I(T^*,I),\\
\frac{\partial E}{\partial t} &=-\nabla \cdot \mathbf{F} - c\sigma_a(T^*)E + ac\sigma_a(T^*)T^4&&\coloneqq \mathbf{N}_E(T^*,E,\mathbf{F},T),\\
    \frac{\partial \mathbf{F}}{\partial t} &=- \frac{c^2}{3}\nabla E - c\sigma_t(T^*)\mathbf F   + \gamma(I,T^*)cE &&\coloneqq \mathbf{N}_F(T^*,I,E,\mathbf{F}),\\
    \frac{\partial T}{\partial t} &= - \frac{1}{\rho c_v}\left(ac\sigma_a(T^*)T^4 - c\sigma_a(T^*)E\right) &&\coloneqq \mathbf{N}_T(T^*,E,T).
\end{align}
\end{subequations}

Note that now we have broken the implicit dependence of $I$ on temperature. As a result, each implicit stage only requires one solve in angular intensity (i.e., sweep), followed by solving the nonlinear LO system. The motivation for this separation of implicit and explicit components is to satisfy the two asymptotic limits of the streaming regime and thick diffusion regime (for stability purposes of IMEX integration), while ensuring the implicit solve only requires one sweep (per stage). In the streaming regime, the transport sweep will resolve the stiff kinetic physics, while in the thick diffusion limit, the nonlinear moment equations correspond to a nonlinear radiation diffusion equation coupled to temperature, which accurately represents the stiff physics. By capturing the two stiff limits of the underlying PDE, we believe that the remaining coupling can be treated explicitly, and will be sufficiently damped by the implicit component to allow stable IMEX/semi-implicit integration. One downside of the SIMEX-RK approach is the double storing of variables, particularly in high-dimensional problems such as transport that already have significant memory requirements. However, note that with our formulation in \eqref{eq:limex-trt1}, we only need $T$ in the $*$-variables, so the additional storage over a standard RK scheme is marginal. It is worth pointing out that if one were to break the implicit coupling in the other direction, solving the LO system first and then performing a transport sweep with the updated temperature profile, one would need to double store all copies of angular flux. Due to its high dimensionality and storage cost, this is a significant downside and we do not pursue here. \tcb{It should be pointed out that within our partitioning, we also treat the opacities semi-implicitly, i.e. linearized about the previous stage. This is a natural extension of the classical semi-implicit integration used in transport simulation \cite{Larsen.1988}. Some semi-implicit schemes and their application to opacities in, e.g. radiation diffusion, are analyzed in \cite{lowrie2004comparison,Ober.2004,knoll2003balanced}, which is related to the much larger class of Rosenbrock methods \cite{rosenbrock1963some}. It is worth pointing out that rather than implicit-explicit, the nonlinear partitioning approach discussed here could also be applied in a semi-implicit form with linearized opacities, facilitating higher-order and adaptive semi-implicit integration. We do not focus on such methods here, as the primary motivation is reducing the number of transport sweeps via implicit-explicit integration.}

In the context of \eqref{eq:limex} and \eqref{eq:limex-trt1}, each stage of our SIMEX-RK formulation of TRT consists of the explicit evaluation, the implicit solve, and the updating of future stages and right-hand sides. Detailed steps for transport-first semi-implicit-explicit integration are presented below.
\begin{enumerate}
    \item \ul{Explicit evaluation:} The first step is to evaluate the explicit part of the operator, $\mathcal{N}_E(Y^*)$ (see Line \ref{line:delta} of \Cref{alg:limex}). In our case, $\boldsymbol{\delta}^I\coloneqq \Delta t\mathcal{N}_E(T^*) = \Delta t\tfrac{1}{4\pi}ac^2\sigma_a(T^*)(T^*)^4$, where $\boldsymbol{\delta}$ is only nonzero in the $I$ equation, and forming $\boldsymbol{\delta}^I$ simply requires a pointwise or element-wise evaluation of the temperature from the $*$ variables.
    
    \item \ul{Implicit solve:} For the $i$th implicit stage \eqref{eq:limex-imp}, assume we have been updating right-hand sides after each implicit solve as in \Cref{alg:limex}, and let $\mathbf{r}_{(i)} \coloneqq \Delta t\sum_{j=1}^{i-1}a_{ij}\mathcal{N}(Y^*_{(j)},Y_{(j)})$ for block solution vector $Y = [I,E,\mathbf{F},T]^T$ and residual $\mathbf{r}_{(i)} = [\mathbf{r}_{(i)}^I, \mathbf{r}_{(i)}^E,\mathbf{r}_{(i)}^{\mathbf{F}},\mathbf{r}_{(i)}^T]^T$. We now do a pre-processing step updating $\mathbf{r}_{(i)}$ with the previous solution value, and $\mathbf{r}^I_{(i)}$ to also include the explicit dependence on $T_{(i)}^*$ (see Line \ref{line:rhs} of \Cref{alg:limex}):
    \begin{equation}\label{eq:Ihat-update}
        \mathbf{r}^I_{(i)} \mathrel{{+}{=}} I_n + a_{ii}\Delta t \tfrac{ac^2}{4\pi}\sigma_a(T_{(i)}^*) T_{(i)}^{4*} \coloneqq I_n + a_{ii}\boldsymbol{\delta}^I.
    \end{equation}
    Then, the implicit set of equations associated with \eqref{eq:limex-trt1} take the block semi-linear form
    {\footnotesize\setlength\arraycolsep{1pt}
    \begin{equation*}
        \begin{bmatrix} I_{(i)} \\ E_{(i)} \\ \mathbf{F}_{(i)} \\ T_{(i)} \end{bmatrix} 
        - a_{ii}\Delta t
        \begin{bmatrix}
        -c\bOmega\cdot\nabla - c\sigma_t(T_{(i)}^*) & \mathbf{0} & \mathbf{0} & \mathbf{0} \\
        \mathbf{0} & -c\sigma_a(T_{(i)}^*) & -\nabla\cdot & ac\sigma_a(T_{(i)}^*)T_{(i)}^3 \\
        \gamma(\cdot ,T_{(i)}^*) & -\tfrac{c^2}{3}\nabla + \gamma c &~ -c\sigma_t(T_{(i)}^*) & \mathbf{0} \\
        \mathbf{0} & \tfrac{c}{\rho c_v}\sigma_a(T_{(i)}^*) & \mathbf{0} & -\tfrac{ac}{\rho c_v}\sigma_a(T_{(i)}^*)T_{(i)}^3
        \end{bmatrix}
        \begin{bmatrix} I_{(i)} \\ E_{(i)} \\ \mathbf{F}_{(i)} \\ T_{(i)} \end{bmatrix}
        = 
        \begin{bmatrix} \mathbf{r}_{(i)}^I \\ \mathbf{r}_{(i)}^E \\ \mathbf{r}_{(i)}^{\mathbf{F}} \\ \mathbf{r}_{(i)}^T \end{bmatrix},
    \end{equation*}
    }where $\gamma(\cdot ,T_{(i)}^*)$ indicates the first function variable is evaluated at the corresponding vector variable $I_{(i)}$.
    
    \item[2a.] \emph{Transport sweep and evaluating $\gamma$:} Note that we can solve the leading equation for $I_{(i)}$ via a transport sweep as it is does not depend on $\{E_{(i)},\mathbf{F}_{(i)}, T_{(i)}\}$, and eliminate it from the system. Following Line \ref{line:N} of \Cref{alg:limex}, redefine 
    \begin{equation}\label{eq:deltaI}
        \boldsymbol{\delta}^I \pluseq \frac{I_{(i)} - \mathbf{r}_{(i)}^I}{a_{ii}},
    \end{equation}
    after which we have $\boldsymbol{\delta}^I = \Delta t\mathbf{N}_I(T_{(i)}^*,I_{(i)})$. From \eqref{eq:limex-trt1} and we then have 
    \begin{align*}
        \frac{\partial I(T_{(i)}^*,I_{(i)})}{\partial t} & = \mathbf{N}_I(T_{(i)}^*,I_{(i)}) 
        \\& = \tfrac{1}{4\pi}ac^2\sigma_a(T_{(i)}^*) T_{(i)}^{4*} - c\Omega \cdot \nabla I_{(i)} - c\sigma_t(T_{(i)}^*) I_{(i)}
        = \frac{1}{\Delta t}\boldsymbol{\delta}^I.
    \end{align*}
    
    We then construct the consistency term $\gamma(I_{(i)},T_{(i)}^*)$ as in \eqref{eq:indirect-partial-F} by taking moments of $\partial I(T_{(i)}^*,I_{(i)})/\partial t = \boldsymbol{\delta}^I/\Delta t$. 

    \item[2b.] \emph{LO solve:} 
    Following the transport sweep, and evaluation of the consistency term, we arrive at the discrete version of the nonlinear LO equations \eqref{eq:moment-eqs3}:
    \begin{subequations}
    \begin{align*}
        \frac{E_{(i)} - \mathbf{r}^E_{(i)}}{a_{ii}\Delta t} & =
            -c\sigma_a(T_{(i)}^*)E_{(i)} -\nabla\cdot\mathbf{F}_{(i)} + ac\sigma_a(T_{(i)}^*)T_{(i)}^4, \\
        \frac{\mathbf{F}_{(i)} - \mathbf{r}^{\mathbf{F}}_{(i)}}{a_{ii}\Delta t} & = 
            -\frac{c^2}{3}\nabla E_{(i)} -c\sigma_t(T_{(i)}^*)\mathbf{F}_{(i)}+\gamma cE_{(i)}, \\
        \frac{T_{(i)} - \mathbf{r}^T_{(i)}}{a_{ii}\Delta t} & = \tfrac{c}{\rho c_v}\sigma_a(T_{(i)}^*)E_{(i)} -\tfrac{ac}{\rho c_v}\sigma_a(T_{(i)}^*)T_{(i)}^4,
    \end{align*}
    \end{subequations}
    where $\gamma$ is fixed for all inner nonlinear LO iterations. This is effectively a backward Euler time step applied to LO moment equations with modified starting solution (residual $\mathbf{r}_{(i)}$) and time step $a_{ii}\Delta t$.

    \item[3.] \ul{Update future stages:} We now complete evaluation of our full nonlinear operator $N(Y_{(i)}^*,Y_{(i)})$ following the implicit solve as in Line \ref{line:N} of \Cref{alg:limex} and \eqref{eq:deltaI} by defining 
    \begin{equation}\label{eq:deltas}
        \boldsymbol{\delta}^E \coloneqq \frac{E_{(i)} - \mathbf{r}_{(i)}^E}{a_{ii}},
        \hspace{3ex}
        \boldsymbol{\delta}^{\mathbf{F}} \pluseq \frac{\mathbf{F}_{(i)} - \mathbf{r}_{(i)}^{\mathbf{F}}}{a_{ii}},
        \hspace{3ex}
        \boldsymbol{\delta}^T \coloneqq \frac{T_{(i)} - \mathbf{r}_{(i)}^T}{a_{ii}},
    \end{equation}
    where $\boldsymbol{\delta}^E = \Delta t\mathbf{N}_E$, $\boldsymbol{\delta}^{\mathbf{F}} = \Delta t\mathbf{N}_F$, and
    $\boldsymbol{\delta}^T = \Delta t\mathbf{N}_T$ from \eqref{eq:limex-trt1}. From here we use $\boldsymbol{\delta}$ to update future stage solutions and right-hand sides as in Line \ref{line:loop} of \Cref{alg:limex}.
\end{enumerate}

\noindent There are a few things to note about the proposed method:
\begin{itemize}
    \item If a domain has reflective boundary conditions, one must be careful when solving the implicit equation. Care must be taken to either (i) perform sweeps in a way such that you exactly invert the transport operator, including boundary conditions, or (ii) break coupling that arises from reflective boundary conditions and treat certain degrees of freedom explicitly, so that a single sweep inverts the implicit transport operator. 
    \item If one wants the robustness of implicit methods and higher than first order accuracy without incorporating opacities or EOS in the implicit solve, the above approach can be modified to be only explicit in opacity, while maintaining implicit coupling between HO and LO systems. This can be seen as a practical generalization of the initial work on DIRK methods for TRT \cite{maginot2016high}\tcb{, or a high-order extension of the first-order semi-implicit approach typically used in transport \cite{lowrie2004comparison,Ober.2004,knoll2003balanced,Larsen.1988}.} We refer to such an approach as ``semi-implicit,'' but do not explore further here, as the computational savings of one sweep is our main objective.
    \item \tcb{There can be interest in (semi-)discrete maximum principles, e.g. \cite{larsen2013properties}. Unlike implicit schemes, our proposed IMEX methods have a finite stability region, and if the timestep is too large and integration unstable, we certainly will not satisfy a maximum principle. Moreover, joint stability of partitioned integration schemes is very difficult even for linear problems, so quantifying the region of stability precisely would be challenging or infeasible. Nevertheless, the LO system is a valid TRT model and integrated (semi-)implicitly, thus it will satisfy a maximum principle without convergence of the consistency term assuming the underlying semi-implicit integrator does. At worst, in say the case where the consistency terms are neglected, the method obeys the discrete diffusion maximum principle. In the case of not converging the consistency term, the scheme obeys a discrete maximum principle that limits to the true maximum principle as the time step and mesh sizes are reduced.}

\end{itemize}

\section{Numerical results}\label{sec:results}
We use the 1d Marshak wave and 2d crooked pipe benchmark problems to demonstrate the convergence and performance of the proposed IMEX-HOLO framework with direct comparison to the traditional, implicit HOLO algorithm from \cite{Park.2020}.  
The Marshak wave problem is a strongly nonlinear problem due to temperature-dependent opacities while the crooked pipe problem is a \tcb{benchmark} multi-material problem that exhibits both advective and diffusive behavior at early and late times, respectively. 
Both problems are stress tests for IMEX stability, typically requiring tuning of parameters and the use of algorithmic hardening functionalities such as negative flux fixups and temperature flooring to run without failing. 
To assess accuracy, each problem is run with a succession of temporal refinements each compared to a time-resolved reference solution. 
All runs use a fixed spatial mesh and angular quadrature order so that temporal error is investigated in isolation. 
The time-resolved reference solution is taken to be the most accurate scheme on that problem using a time step $8\times$ smaller than the smallest time step size used to assess convergence. 
This procedure allows investigating time integration error on challenging problems where an exact solution may not be known. 
Unless otherwise noted, the error is computed using a \emph{common} reference method so that each method converging to the reference indicates the methods are all converging to the true, time-resolved solution for that given mesh and angular quadrature rule. 

We compare the above schemes to the implicit HOLO algorithm from \cite{Park.2020}. 
This algorithm uses a backward Euler time integration scheme for the high and low-order variables as well as the coupling between them and thus requires iterating between high and low-order. 
This leads to a nested iteration where the outer iteration converges the coupling between high and low-order and the inner iteration converges the nonlinear coupling between the closed low-order system and the temperature equation.  
\tcb{As in IMEX-HOLO, the temperature is nonlinearly resolved with a pointwise Newton iteration.}
In this study, we have seen that this typically incurs the cost of 5-30 sweeps per time step and 3-7 low-order, linear solves per sweep, meaning that IMEX-HOLO's one sweep and one nonlinear low-order solve per stage will significantly reduce the total cost of the simulation. 
Following \cite{Park.2020}, the implicit HOLO algorithm fixes the opacities at the beginning of each time step.  
Note that our SIMEX schemes generalize this such that opacities are fixed at each stage. 
Finally, the implicit algorithm uses the finite difference approximation for the time derivative in the HOLO consistency term as described in \Cref{sec:holo:gamma}. 

All simulations are run on compute nodes with dual socket Xeon Gold 6152 22-core processors.

\subsection{1d Marshak problem}\label{sec:results:marshak}
\begin{figure}
\centering
\includegraphics[width=.5\textwidth]{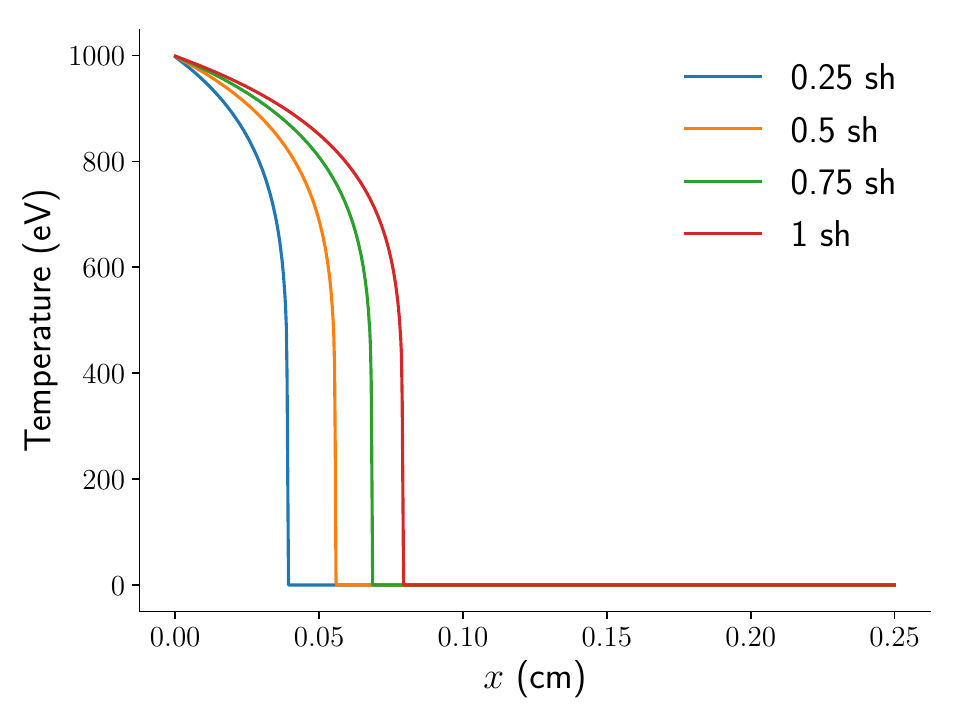}
\caption{Evolution of the temperature on the Marshak wave problem over time. Radiation impinges on the slab from the left. As the initially cold material heats, it becomes transparent, allowing radiation to penetrate further into the domain. The steep temperature gradient at the front of the wave poses significant challenges for robustness. Note that $\SI{1}{\sh} = \SI{10}{\nano\second}$.}
\label{fig:marshak_evolution}
\end{figure}

The Marshak wave problem models radiation impinging on a slab where the slab's material properties are strongly temperature dependent. 
The computational domain is $[0,\SI{0.25}{\cm}]$ and the simulation time is $t = \SI{10}{\ns}$. 
The temperature and radiation fields are initially in equilibrium at a temperature of \SI{.025}{\eV}. 
At time $t=0$, a strong radiation source at $x=0$ corresponding to a Planckian distribution at $T = \SI{1}{\keV}$ is turned on. 
The material opacity is $\sigma(T) = 10^{12}/T^3\si{\per\cm}$ and the heat capacity is $C_v = \SI{3e12}{\erg\per\eV\per\cm\cubed}$. 
Initially, the material is cold and thus extremely optically thick with $\sigma(\SI{0.025}{\eV}) = \mathcal{O}(10^{16})$.
As the material heats, the material becomes transparent to radiation allowing the radiation to penetrate further into the domain and heat more of the material.  
This creates a wave with a strong temperature gradient that moves through the domain over time. 
Solution profiles at a range of snapshots in time are shown in \Cref{fig:marshak_evolution}. 
We use a mesh of 1000 cells and S$_8$ angular quadrature. 

\begin{figure}
\centering
\begin{subfigure}{.32\textwidth}
    \centering
    \includegraphics[width=\textwidth]{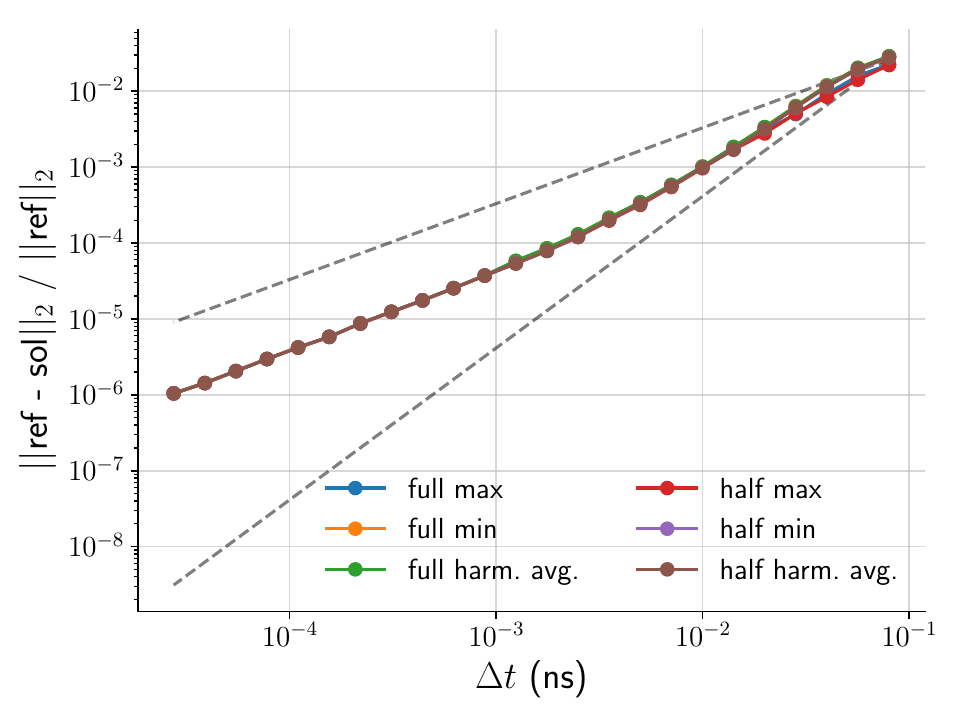}
    \caption{}
\end{subfigure}
\begin{subfigure}{.32\textwidth}
    \centering
    \includegraphics[width=\textwidth]{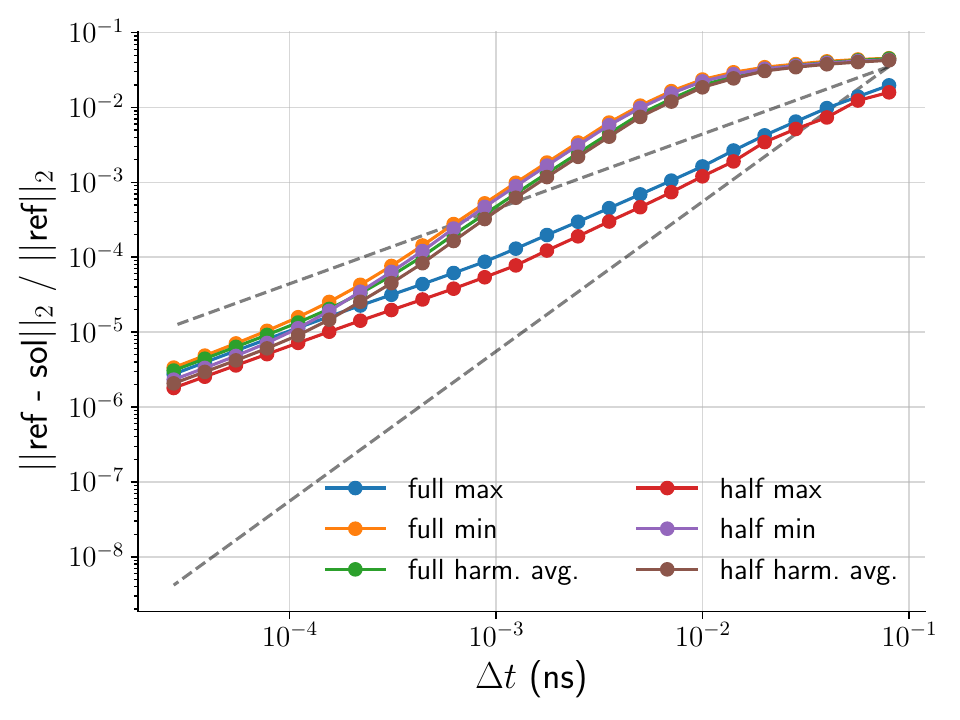}
    \caption{}
\end{subfigure}
\begin{subfigure}{.32\textwidth}
    \centering
    \includegraphics[width=\textwidth]{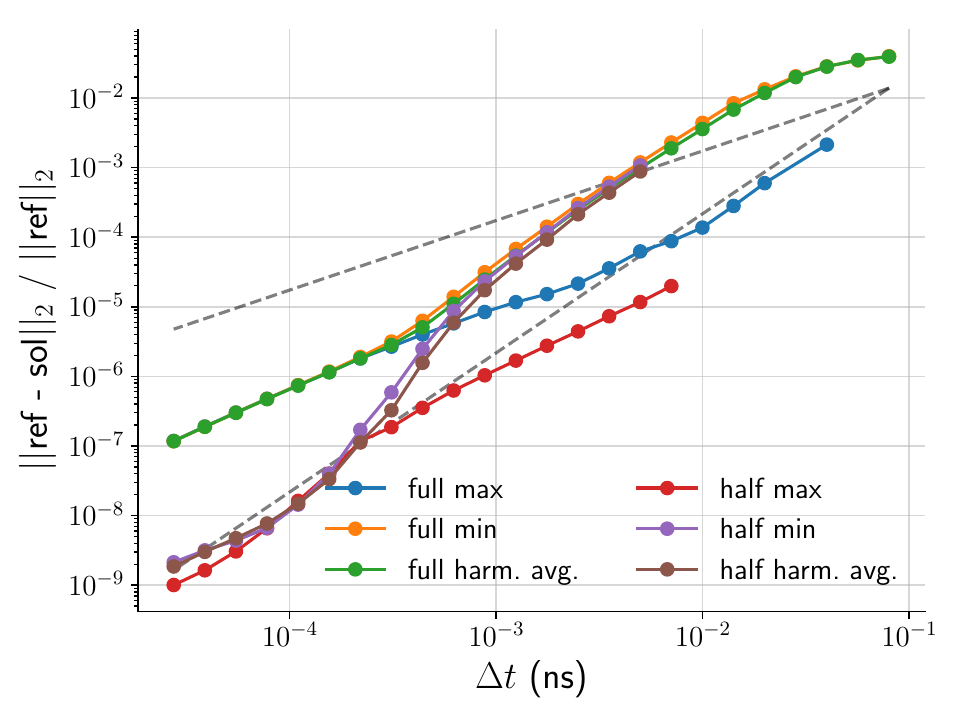}
    \caption{}
\end{subfigure}
\caption{Convergence to a reference solution for the radiation energy density as a function of the time step size, $\Delta t$, for (a) implicit, (b) LIMEX-Euler, and (c) H-LDIRK2(2,2,2). A range of algorithmic choices are compared corresponding to the use of full or half range boundary conditions and the maximum, minimum, and harmonic average value for the interface opacity. \tcb{Gray dashed lines denote first- and second-order convergence. Missing data points indicate a method was not stable at that time step size.}}
\label{fig:compare_bcs_opac_marshak}
\end{figure}
We first compare stability and accuracy with respect to a reference solution for a range of HOLO-related algorithmic choices. 
Here, we use the time-resolved solution generated by H-LDIRK2(2,2,2) as the reference. 
A separate reference is used for each combination of boundary condition type and interface opacity type to avoid biasing the error toward a certain set of parameters. 
\Cref{fig:compare_bcs_opac_marshak} shows the error with respect to the reference for the implicit, LIMEX-Euler, and H-LDIRK2(2,2,2) schemes for all the combinations of boundary condition and interface opacity treatment. 
Implicit converges the coupling between high and low-order, resulting in high and low-order solutions that are consistent to the iteration tolerance. 
In this way, implicit produces the same solution, up to iterative tolerances, regardless of boundary conditions or interface opacity treatment. 
In IMEX-HOLO, the high and low-order systems are no longer temporally consistent due to the explicit and implicit evaluations of the emission term in the high and low-order systems, respectively; the high and low-order systems are equivalent only in the limit as $\Delta t \rightarrow 0$. 
For LIMEX-Euler with moderate time step sizes, the maximum interface opacity option leads to $10\times$ more accuracy compared to the minimum or harmonic average option for both boundary condition types. 
This effect is reduced on the smallest time step sizes. 
H-LDIRK2(2,2,2) is even more sensitive, with the full range boundary conditions resulting in order reduction to first-order accurate. 
The maximum interface opacity option again leads to increased accuracy at moderate time step sizes with minimum and harmonic average lagging up to two orders of magnitude behind. 
While the most accurate, reduced stability is observed with half range boundary conditions or the maximum interface opacity. 

\begin{figure}
\centering
\begin{subfigure}{.49\textwidth}
\centering
\includegraphics[width=\textwidth]{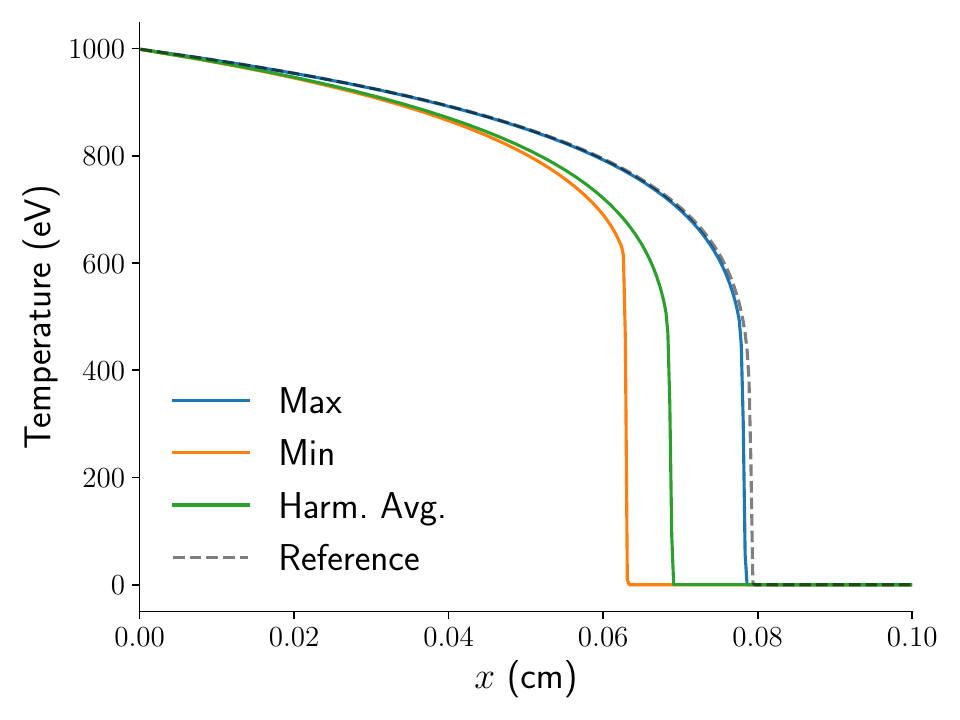}
\caption{}
\end{subfigure}
\begin{subfigure}{.49\textwidth}
\centering
\includegraphics[width=\textwidth]{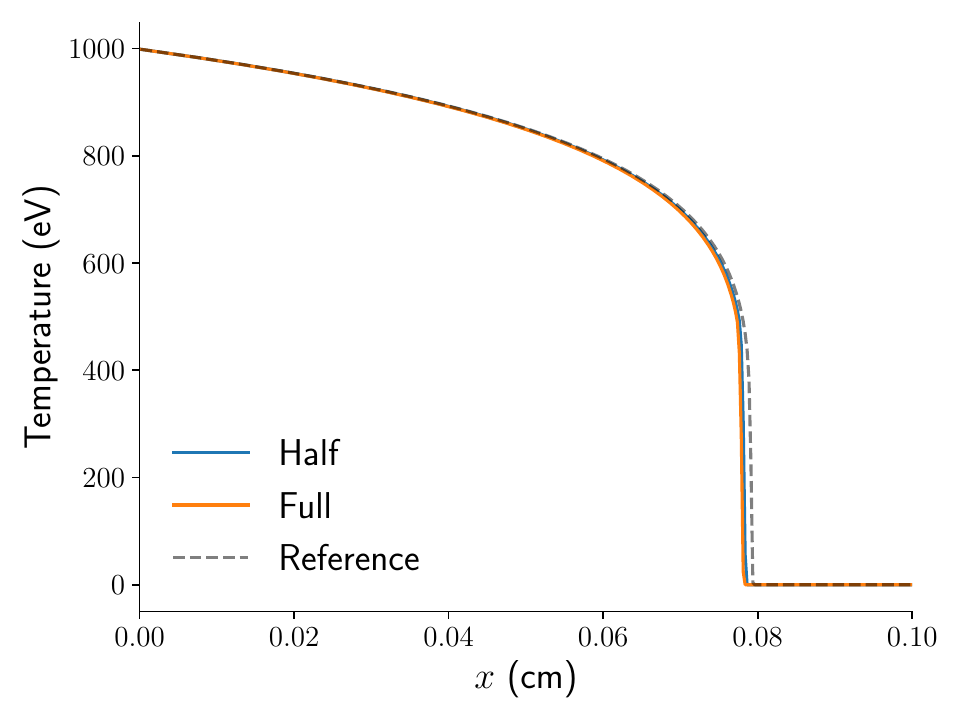}
\caption{}
\label{fig:marshak_sol_compare_bcs}
\end{subfigure}
\caption{The final temperature profile at $t=\SI{10}{\ns}$. In (a), the interface opacity options are varied with the boundary condition type fixed to half range while (b) fixes the interface opacity type to maximum and varies the boundary condition type. All solutions used a moderate time step of size \SI{8e-3}{\ns}. Implicit HOLO with a time step of size \SI{3e-5}{\ns} is used as the reference solution.}
\label{fig:marshak_sol_compare}
\end{figure}

Solution profiles at time $t=\SI{10}{\ns}$ generated with a moderate time step of size \SI{8e-3}{\ns} are compared in \Cref{fig:marshak_sol_compare}. 
The implicit HOLO method with a time step of size \SI{3e-5}{\ns} is used as a reference. 
With the choice of boundary condition fixed at half range, the interface opacity options change both the shape and location of the wave. 
Compared to the reference, the maximum, minimum, and harmonic average options produce relative errors in the final temperature solution of 5\%, 35\%, and 26\%, respectively. 
The maximum option results in a solution closest to the reference.
Choosing the maximum opacity at the interface between cold and hot selects the larger, cold opacity value causing more energy to be deposited into the material. 
This in turn results in the material heating up earlier and radiation penetrating through the domain more rapidly. 
Choosing the minimum opacity selects the smaller, hot opacity, resulting in less energy passing through the interface and an overall slower wave. 
Harmonic average splits the difference between these behaviors. 
An analogous comparison is shown in \Cref{fig:marshak_sol_compare_bcs} where the boundary conditions are varied with a fixed choice of the maximum interface opacity. 
Here, the waves are almost identical having errors with respect to the reference of 5\% and 6\% for half and full range boundary conditions, respectively. 
It may be the case that the full range boundary conditions are in some sense stiffer than the half range boundary conditions resulting in order-reduction to H-LDIRK2(2,2,2)'s first-order stage order. 

\begin{figure}
\centering
\begin{subfigure}{.49\textwidth}
    \centering
    \includegraphics[width=\textwidth]{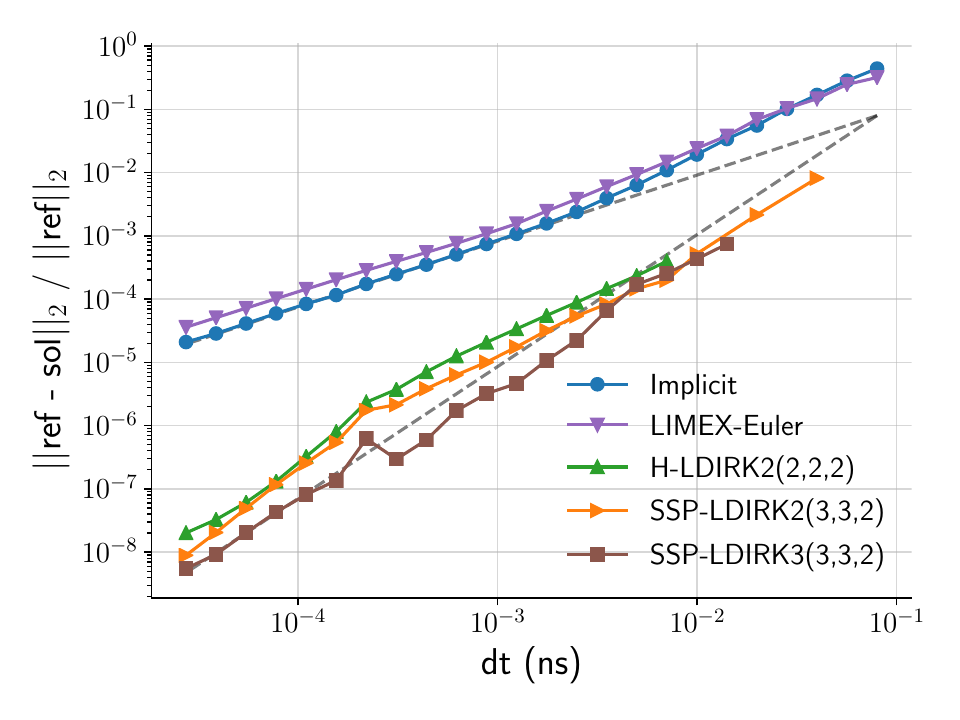}
    \caption{}
\end{subfigure}
\begin{subfigure}{.49\textwidth}
    \centering
    \includegraphics[width=\textwidth]{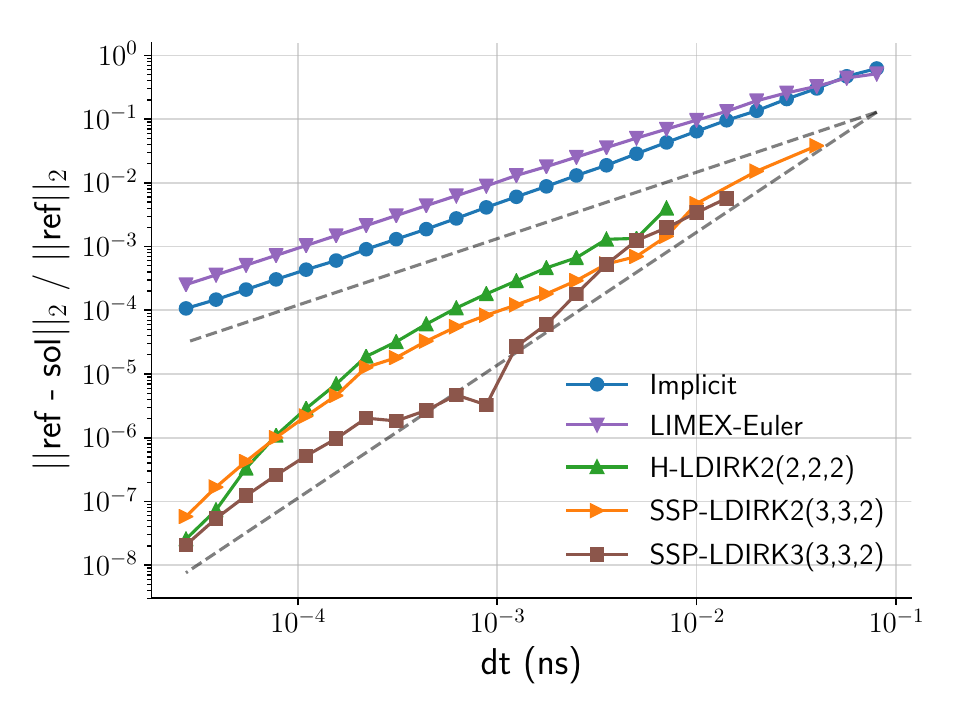}
    \caption{}
\end{subfigure}
\caption{Error with respect to the reference solution as the time step is decreased in the (a) radiation energy density and (b) temperature for a range of HOLO schemes. The time-resolved solution generated by SSP-LDIRK3(3,3,2) is used as the reference. \tcb{Gray dashed lines denote first- and second-order convergence. Missing data points indicate a method was not stable at that time step size.}}
\label{fig:marshak_convergence}
\end{figure}

In the results that follow, we elect to use half range boundary conditions and the maximum interface opacity as this combination was the most accurate for both LIMEX-Euler and H-LDIRK2(2,2,2). 
It is possible that this choice is problem dependent and may require tuning if applied to other problems. However, new integration schemes we are developing \cite{nprk1} applied to half range boundary conditions and the maximum interface opacity offer improved stability, approaching that of LIMEX-Euler, with comparable 2nd-order accuracy. Using these parameters, we now directly compare the accuracy of the implicit, LIMEX-Euler, H-LDIRK2(2,2,2), SSP-LDIRK2(3,3,2), and SSP-LDIRK3(3,3,2) schemes. 
Convergence to the reference is shown in \Cref{fig:marshak_convergence}. 
All methods show optimal convergence with implicit and LIMEX-Euler converging at first-order and H-LDIRK2(2,2,2), SSP-LDIRK2(3,3,2), and SSP-LDIRK3(3,3,2) converging at second-order, as expected. 
LIMEX-Euler has only a minor reduction in accuracy compared to implicit while only performing one sweep per time step. 
In addition, the second-order schemes are significantly more accurate than implicit. 
SSP-LDIRK2(3,3,2), a three-stage, second-order scheme, is the most stable of the high-order schemes, followed by the three stage, second-order SSP-LDIRK3(3,3,2), and finally the two-stage, second-order H-LDIRK2(2,2,2). 
This indicates that methods with additional stages can offer both improved stability and accuracy. 
SSP-LDIRK3(3,3,2) is the most accurate, possibly due to its use of a third-order explicit scheme to evaluate $T^*$, leading to more accurate opacities and emission sources than H-LDIRK2(2,2,2) and SSP-LDIRK2(3,3,2) which use second-order explicit schemes for $T^*$. This scheme was also consistently the most accurate in our gray radiation hydrodynamics problems as well \cite{southworth2023implicit}.

\begin{figure}
\centering
\begin{subfigure}{.32\textwidth}
    \centering
    \includegraphics[width=\textwidth]{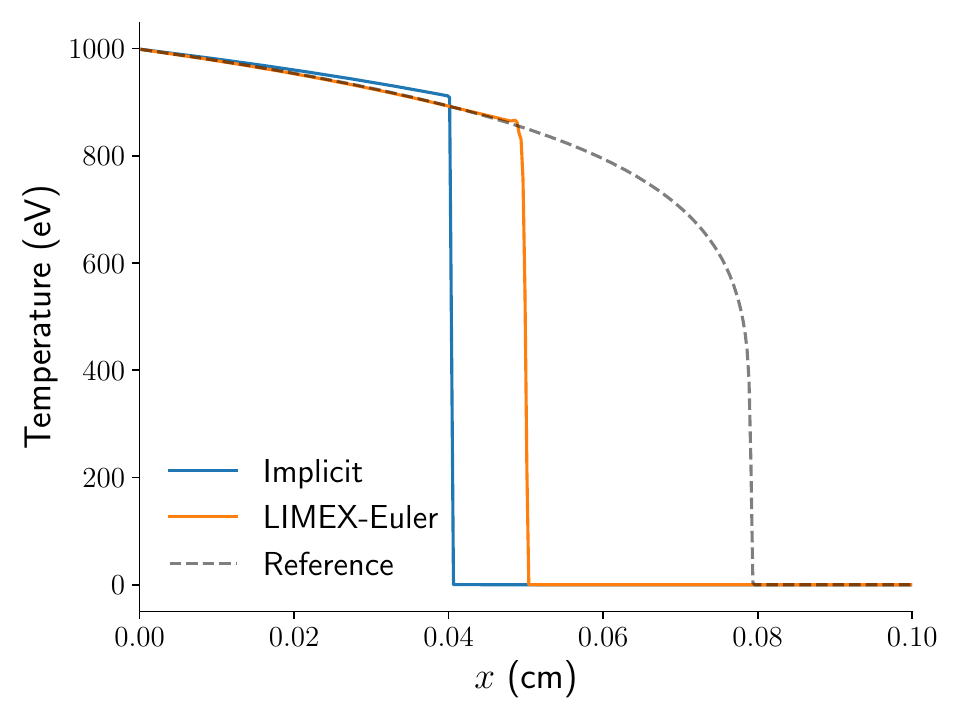}
    \caption{}
\end{subfigure}
\begin{subfigure}{.32\textwidth}
    \centering
    \includegraphics[width=\textwidth]{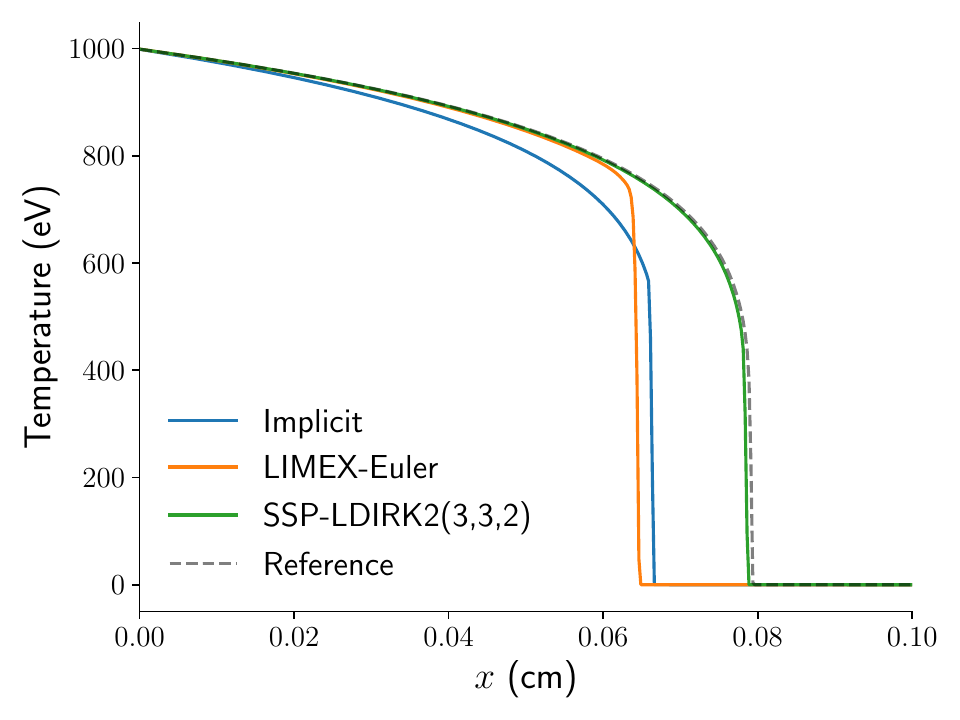}
    \caption{}
\end{subfigure}
\begin{subfigure}{.32\textwidth}
    \centering
    \includegraphics[width=\textwidth]{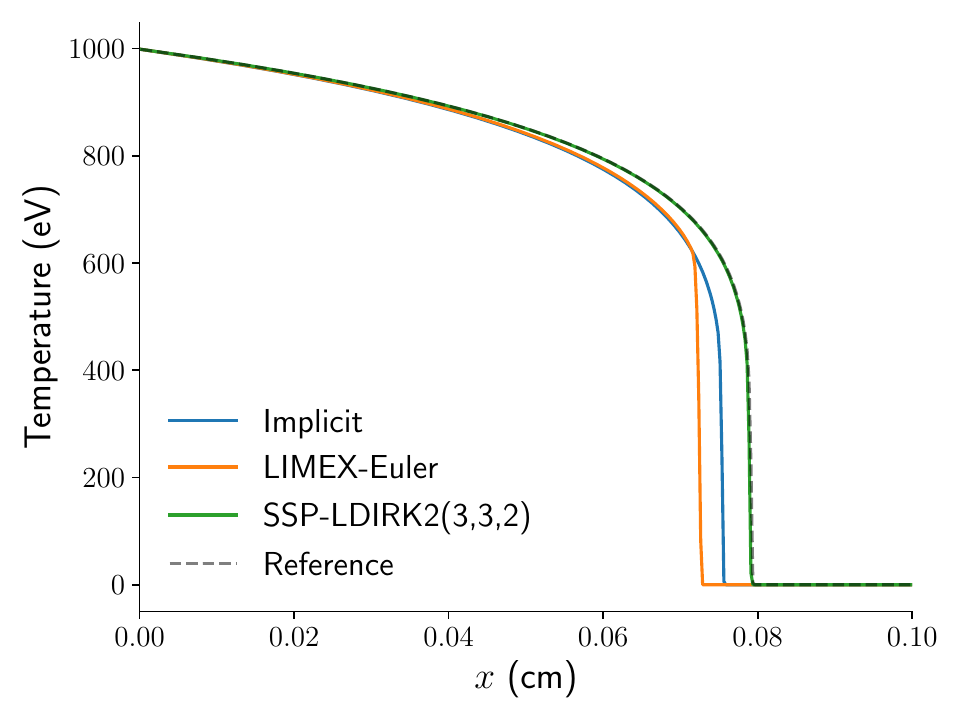}
    \caption{}
\end{subfigure}
\caption{Temperature profiles at $t=\SI{10}{\ns}$ generated with time steps of size (a) \SI{8e-2}{\ns}, (b) \SI{4e-2}{\ns}, and (c) \SI{2e-2}{\ns}.}
\label{fig:marshak_final_waves}
\end{figure}

Finally, we compare solution quality when large time steps are used. 
\Cref{fig:marshak_final_waves} shows the final temperature profile generated by simulations with time steps of size \SI{8e-2}{\ns}, \SI{4e-2}{\ns}, and \SI{2e-2}{\ns}. 
Note that only some of the methods are stable with such large time steps sizes and are thus not plotted in \Cref{fig:marshak_final_waves}. 
Observe that even with one sweep per time step, LIMEX-Euler produces a temperature profile comparable with implicit. 
Additionally, when stable, SSP-LDIRK2(3,3,2) is visually producing the reference solution with three sweeps per time step, even with relatively large time steps.

\subsection{2d crooked pipe problem}\label{sec:results:tophat}
\begin{figure}[thb!]
    \centering
    \includegraphics[width=0.75\textwidth]{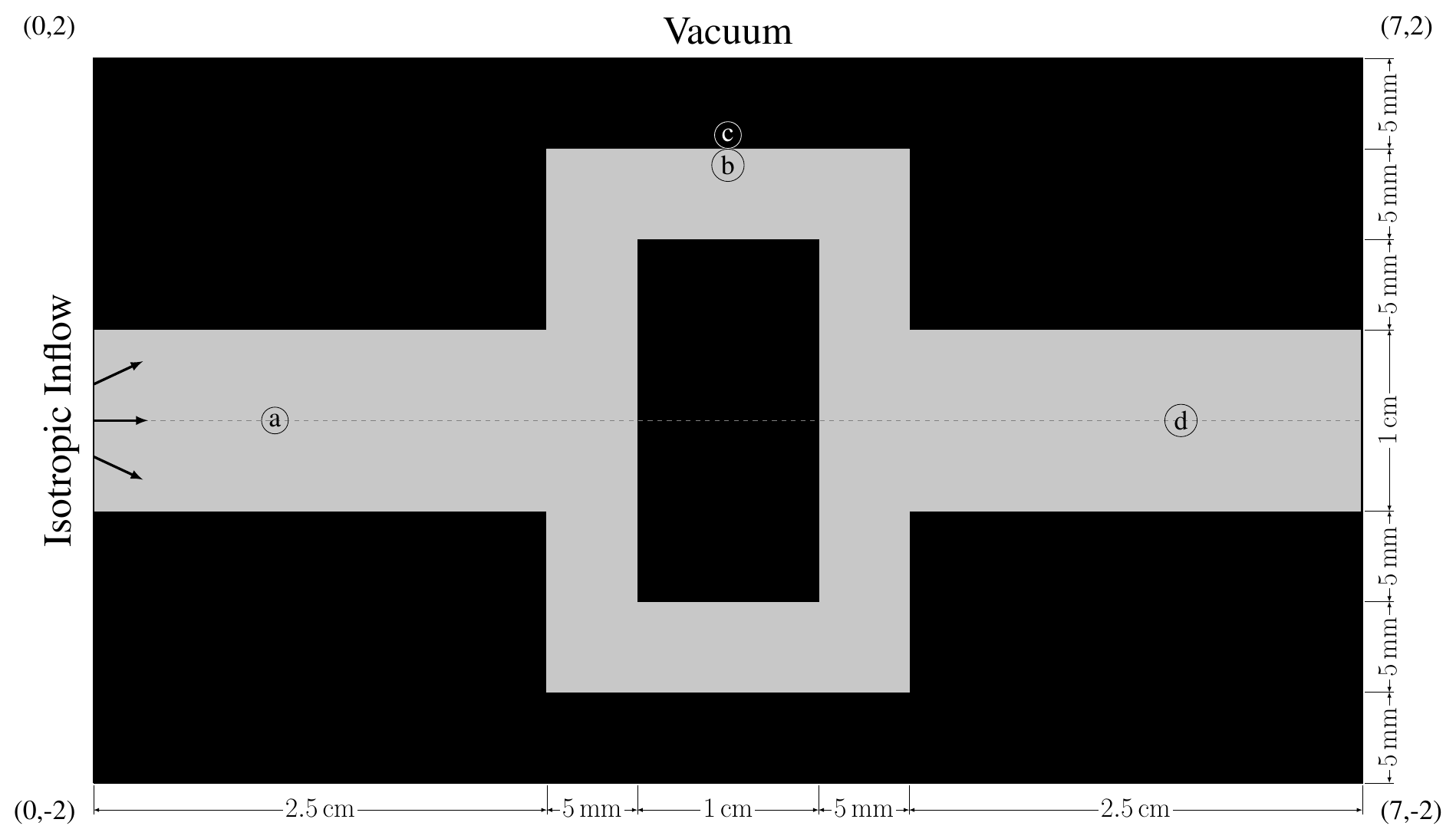}
    \caption{Geometry and materials for the crooked pipe benchmark. This benchmark consists of two materials: the optically thin pipe (gray) and optically thick wall (black). Radiation isotropically enters the pipe on the left side of the domain. A reflection plane at $y=0$ is used to halve the computational domain with no inflow, vacuum boundary conditions applied on the remainder of the domain's boundary. Wave arrival speed is monitored with high resolution in time at the four spatial points labeled in the diagram with circled letters. }
    \label{fig:tophat_geometry}
\end{figure}
\begin{figure}
    \centering
    \def\figheight{3in}
    \includegraphics[height=\figheight, trim={0 0 872 65}, clip]{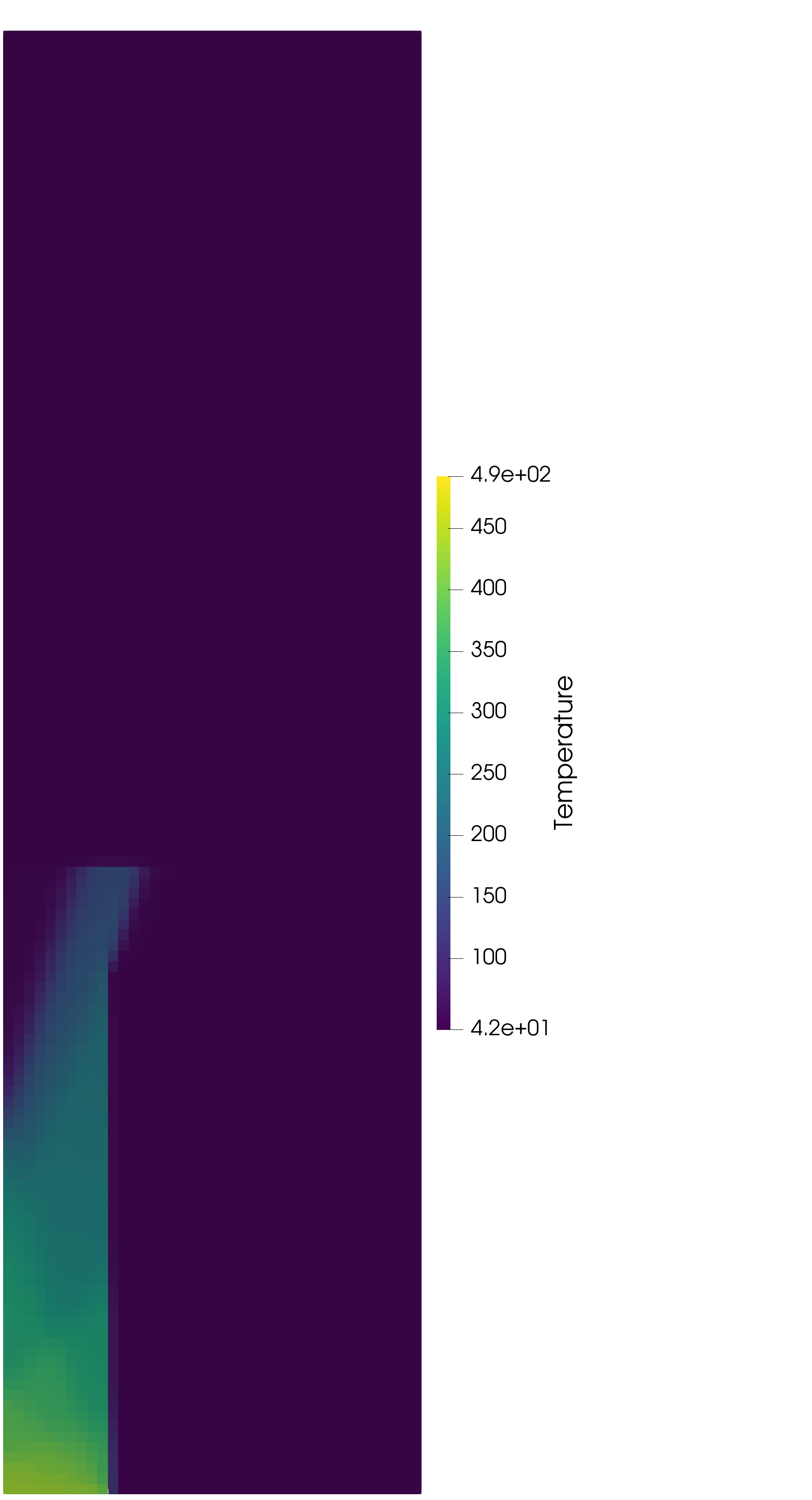}
    \includegraphics[height=\figheight, trim={0 0 872 65}, clip]{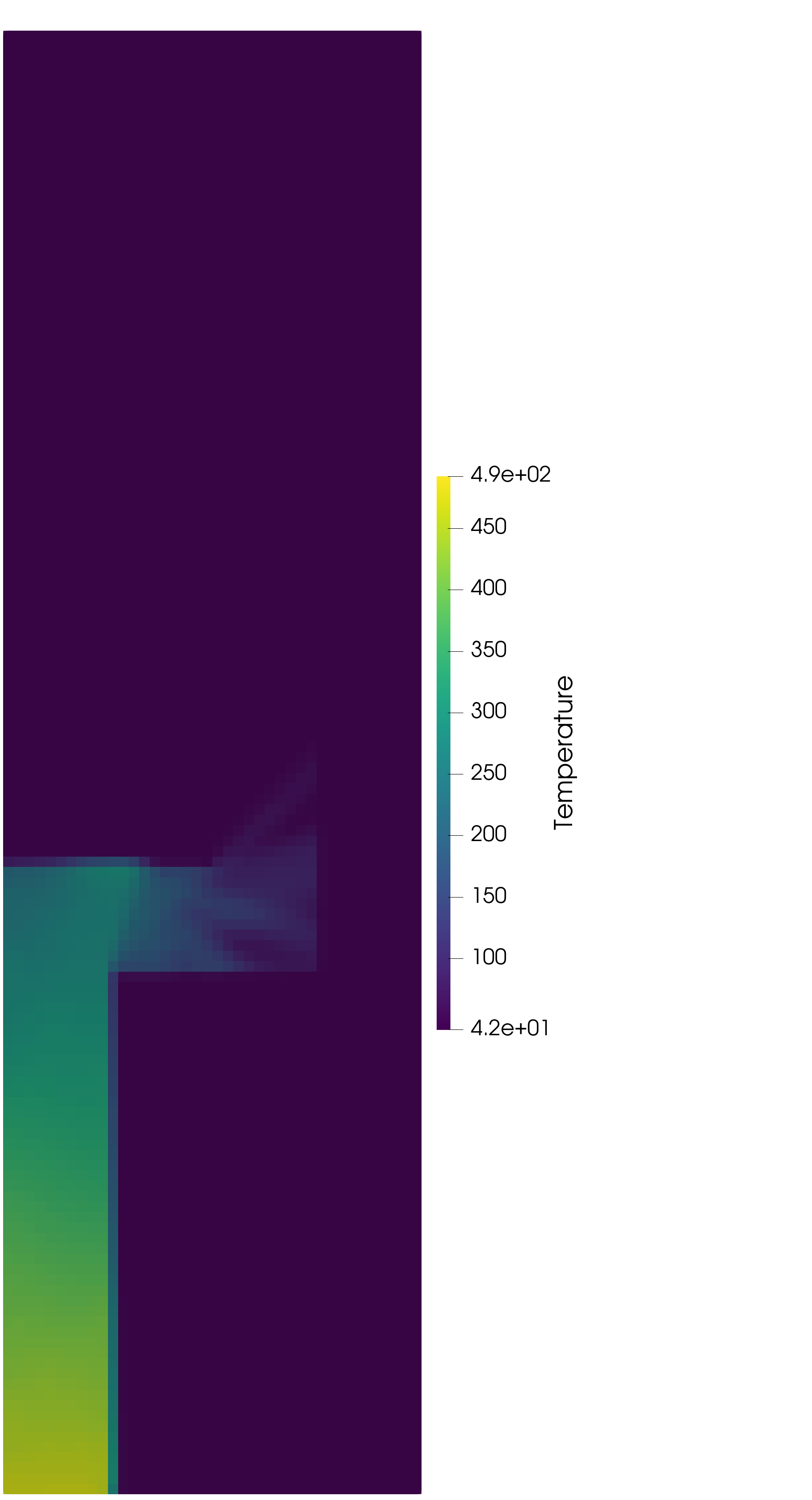}
    \includegraphics[height=\figheight, trim={0 0 872 65}, clip]{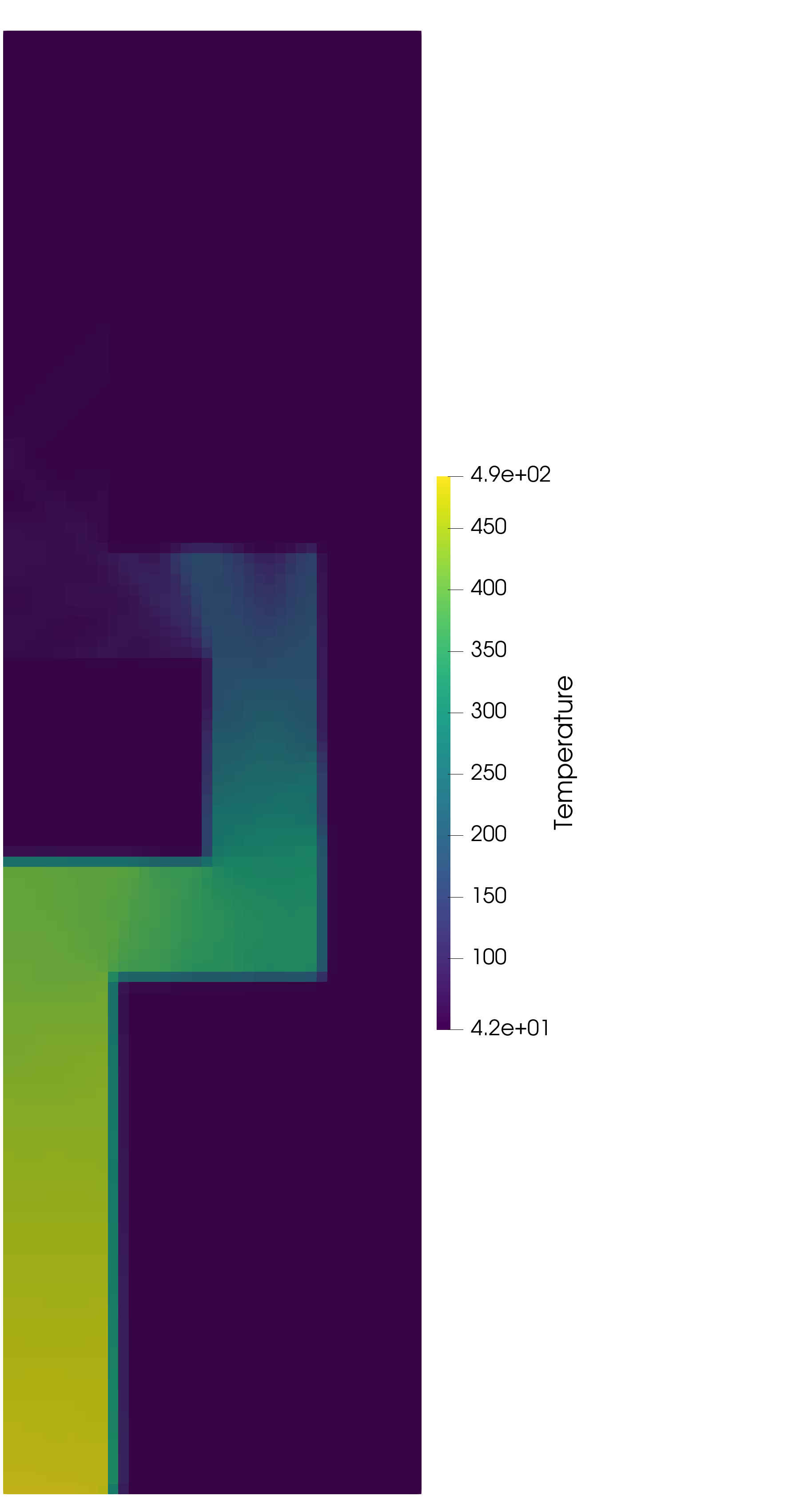}
    \includegraphics[height=\figheight, trim={0 0 872 65}, clip]{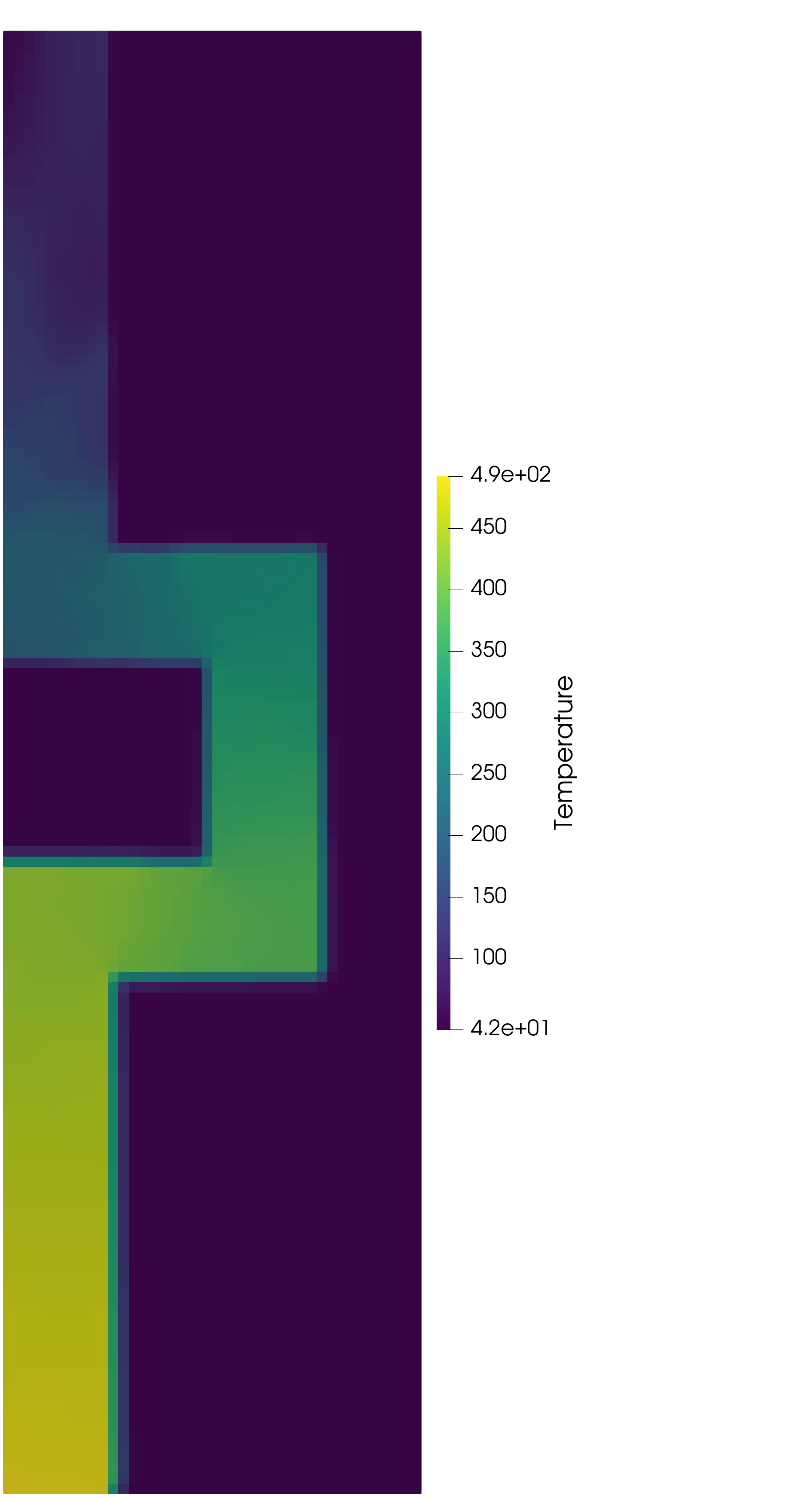}
    \includegraphics[height=\figheight, trim={0 0 525 65}, clip]{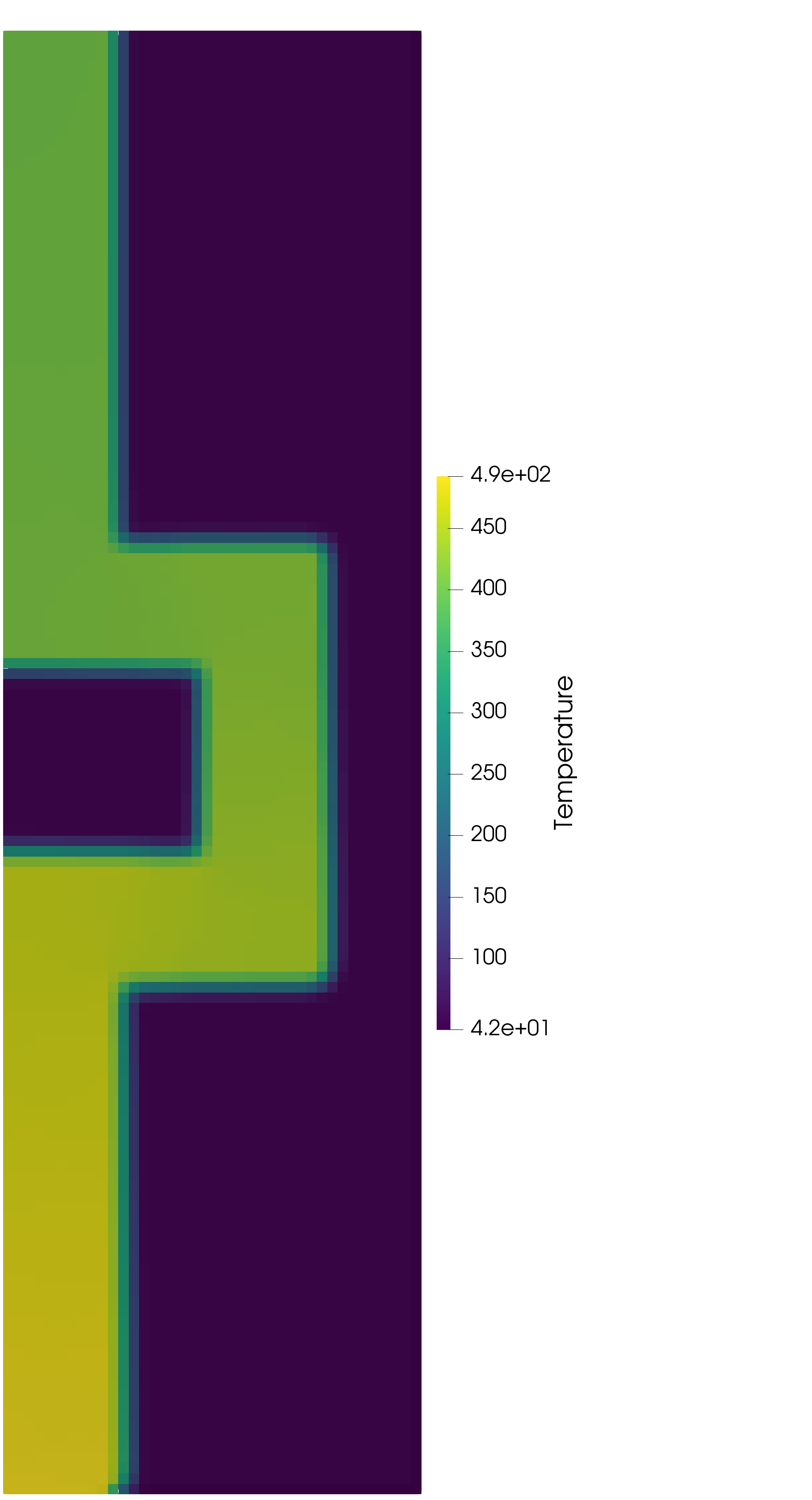}
    \caption{The evolution of the temperature at five snapshots in time. From left to right, the profiles correspond to simulation times of \SI{0.1}{\sh}, \SI{1}{\sh}, \SI{5}{\sh}, \SI{10}{\sh}, and \SI{50}{\sh}, respectively. Level symmetric S$_6$ angular quadrature was used with a uniform spatial mesh of cells of size $\SI{.05}{\cm} \times \SI{.05}{\cm}$.}
    \label{fig:tophat_temperature_evolution}
\end{figure}

We now demonstrate accuracy and computational performance on the multi-material, two-dimensional crooked pipe benchmark problem. 
The geometry and materials are depicted in \Cref{fig:tophat_geometry}. 
The problem consists of two materials: the optically thin pipe and optically thick wall described by 
    \begin{equation}
        \sigma_\text{pipe} = \SI{0.2}{\per\cm} \,, \quad \sigma_\text{wall} = \SI{2000}{\per\cm} \,. 
    \end{equation}
Note that, unlike the Marshak problems, the opacities do not depend on temperature and are thus fixed in time. 
The heat capacity is $C_v = 10^{12}\si{\erg\per\eV\per\cm\cubed}$ in both the wall and pipe. 
Radiation enters the pipe at the left edge of the domain according to a Planckian distribution evaluated at a temperature of \SI{500}{\eV}. 
Elsewhere, the domain boundary is treated as a vacuum. 
In other words, 
    \begin{equation}
        I_\text{inflow} = \begin{cases}
            \frac{ac T_\text{inflow}^4}{4\pi} \,, & x=0, y \in [\SI{-2}{\cm}, \SI{2}{\cm}] \\ 
            0 \,, & \text{otherwise} 
        \end{cases} \,, 
    \end{equation}
with $T_\text{inflow} = \SI{500}{\eV}$. 
The computational domain is halved by applying a reflection boundary along the plane $y=0$. 
The initial temperature and radiation fields are set to be in equilibrium with each other at a spatially uniform temperature of \SI{50}{\eV}. 
\Cref{fig:tophat_geometry} also labels four spatial locations where the pointwise solution is monitored with high resolution in time (i.e.~at every time step) to assess wave arrival time at the beginning, middle, and end of the pipe. 
Note that (b) and (c) are placed at the center of the cells that straddle the pipe-wall interface such that (b) is just inside the pipe and (c) just outside. 
An example of the numerical evolution of the temperature field on this problem is provided in \Cref{fig:tophat_temperature_evolution}. 
Observe that the wave begins to ``turn the corner'' around \SI{0.1}{\sh} and that the solution begins to behave diffusively around \SI{10}{\sh}. 
Thus, we target a maximum time step size of \SI{0.1}{\sh} so that the early-time transport dynamics are always adequately resolved and choose a final simulation time of \SI{50}{\sh} to be able to assess the accuracy and computational performance of the numerical schemes after the problem has become diffusive. 
For all results in this section, S$_6$ Level Symmetric angular quadrature is used with a uniform spatial mesh of cells of size $\SI{.05}{\cm} \times \SI{.05}{\cm}$. 
We use a range of practical and realistic time step sizes between \SI{0.1}{\sh} and $10^{-3}\,\si{\sh}$. 
Given a characteristic spatial mesh length of \SI{.05}{\cm} and a speed of light of \SI{300}{\cm\per\sh}, these time step sizes correspond to advective CFL conditions between 600 and 9.375, meaning that an explicit treatment of the streaming and collision operator would not be stable. 
Note that the reference solution's time step has an advective CFL of $\approx\!0.6$ and thus does resolve the dynamics of advection at the speed of light, although explicit integration would still be unstable for the full transport problem. 

\subsubsection{Convergence to Reference Solution}
\begin{figure}
    \centering
    \begin{subfigure}{.49\textwidth}
        \centering
        \includegraphics[width=\textwidth]{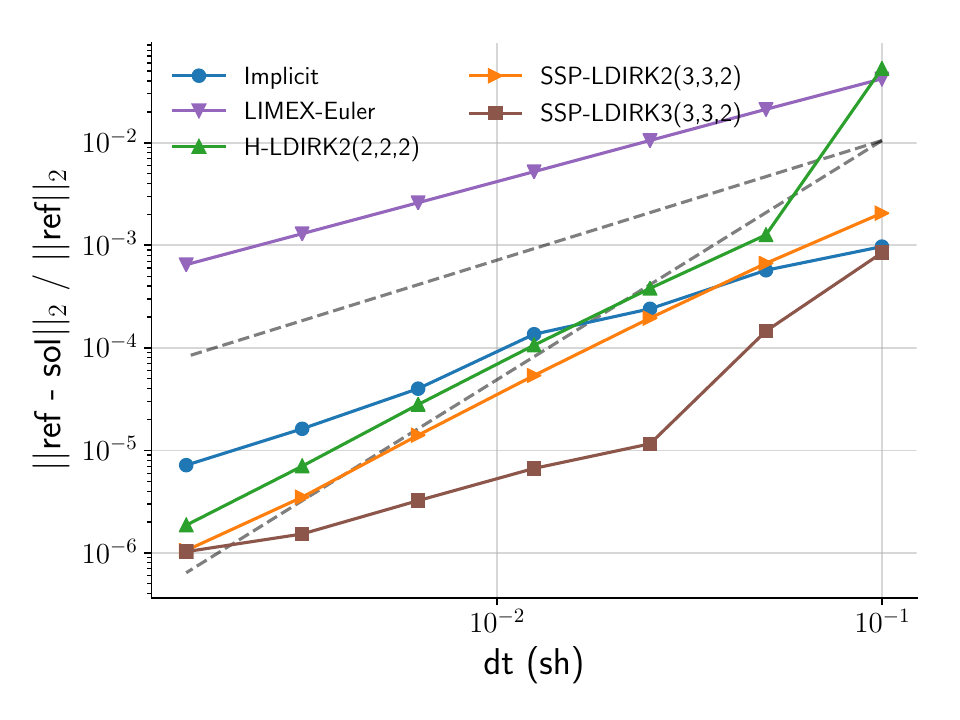}
        \caption{}
    \end{subfigure}
    \begin{subfigure}{.49\textwidth}
        \centering
        \includegraphics[width=\textwidth]{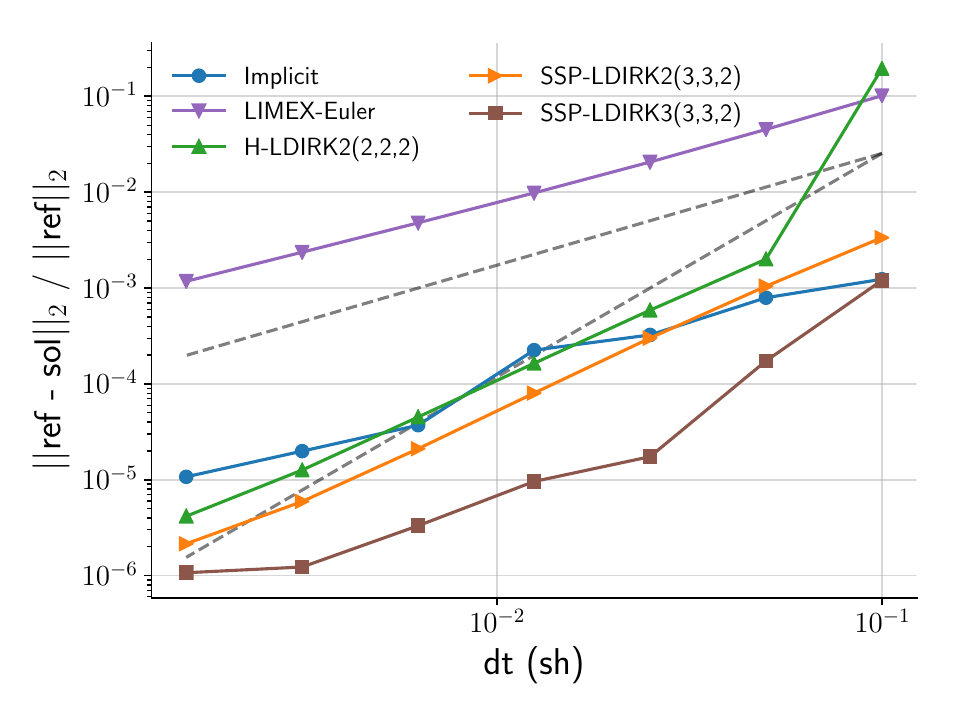}
        \caption{}
    \end{subfigure}
    \caption{Error for integration schemes on the crooked pipe problem for (a) temperature and (b) radiation energy. \tcb{Gray dashed lines denote first- and second-order convergence.}}
    \label{fig:tophat_convergence}
\end{figure}

\Cref{fig:tophat_convergence} shows convergence of the temperature and radiation energy for each scheme to an SSP-LDIRK3(3,3,2) reference solution. 
Using logarithmic regression on all seven radiation energy density error data points, implicit and LIMEX-Euler have slopes of 1.16 and 1.0, respectively. 
H-LDIRK2(2,2,2) struggles in accuracy at the largest time step, corresponding to a single time step before radiation hits the first corner in the domain, showing an anomalously large error for that time step only. 
We note that this is an artifact of prioritizing accuracy over stability in choosing the boundary condition and interface opacity treatment as only the maximum interface opacity treatment led to the observed jump in error for $\Delta t = \SI{0.1}{\sh}$. 
In addition, SSP-LDIRK3(3,3,2) stalls in accuracy around $10^{-6}$. 
This is likely due to the time integration error becoming smaller than other errors present in the simulation, such as errors associated with iterative tolerances. 
Excluding the first point for H-LDIRK2(2,2,2) and the last point for SSP-LDIRK3(3,3,2), H-LDIRK2(2,2,2) and SSP-LDIRK3(3,3,2) achieve orders of 1.8 and 1.93, respectively. 
Thus, all methods are converging near their expected, optimal rates. 

The IMEX schemes have larger error constants relative to implicit than seen on the Marshak problem. 
Where LIMEX-Euler and implicit had comparable error on the Marshak problem, LIMEX-Euler is now two orders of magnitude less accurate than the implicit scheme. 
The second-order H-LDIRK2(2,2,2) scheme is able to beat the implicit method in accuracy only for time steps smaller than \SI{.015}{\sh} with comparable accuracy otherwise aside from the anomalous first data point. 
On the crooked pipe, SSP-LDIRK3(3,3,2)'s extra stage has a more pronounced improvement in accuracy, often beating H-LDIRK2(2,2,2) by over an order of magnitude, similar to our experience in radiation hydrodynamics \cite{southworth2023implicit}.
\begin{figure}
    \centering
    \begin{subfigure}{.49\textwidth}
        \centering
        \includegraphics[width=\textwidth]{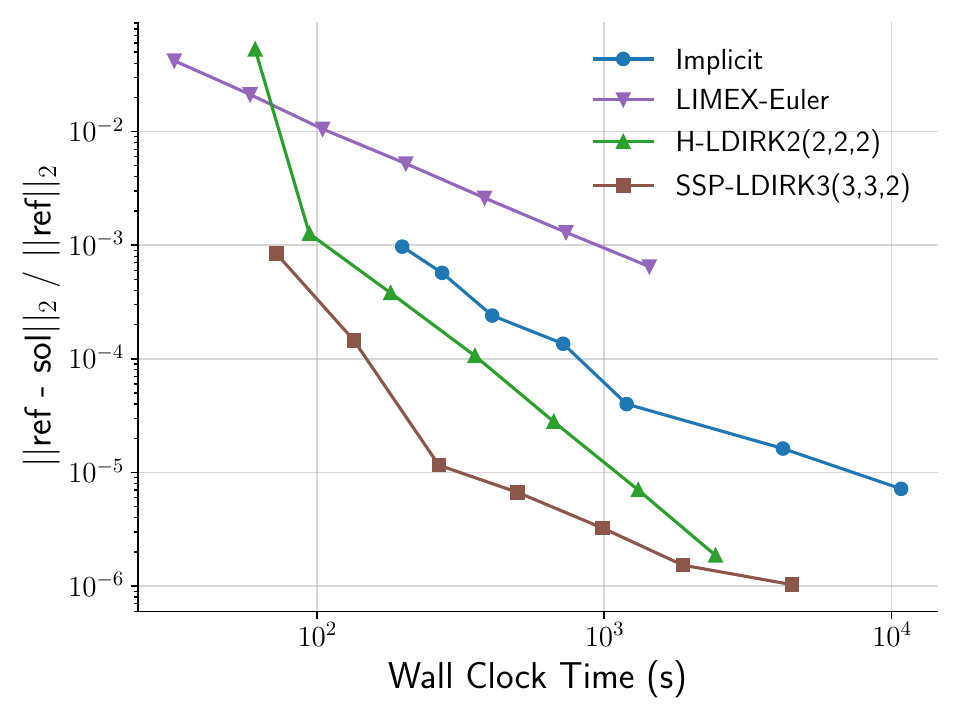}
        \caption{}
    \end{subfigure}
    \begin{subfigure}{.49\textwidth}
        \centering
        \includegraphics[width=\textwidth]{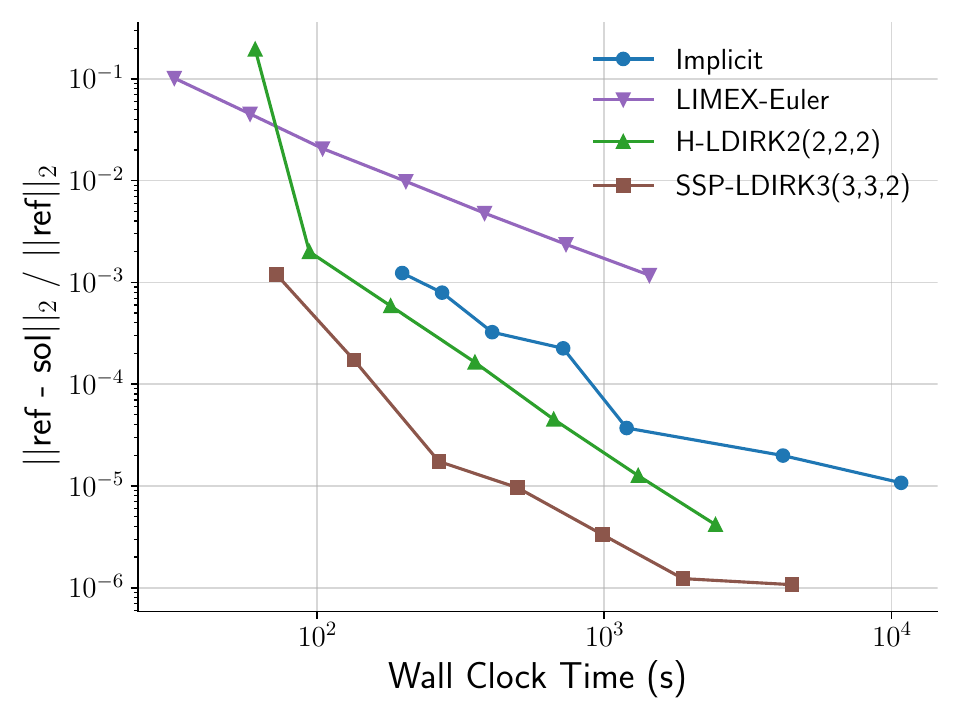}
        \caption{}
    \end{subfigure}
    \caption{Accuracy in computing (a) temperature and (b) radiation energy as a function of total simulation wall clock time for IMEX shemes compared to implicit.}
    \label{fig:tophat_efficiency}
\end{figure}
Accuracy as a function of total simulation wall clock time is presented in \Cref{fig:tophat_efficiency}. 
For a given accuracy, SSP-LDIRK3(3,3,2) has the lowest time-to-solution, followed by H-LDIRK2(2,2,2), implicit, and finally LIMEX-Euler. This indicates that, on this problem, LIMEX-Euler's one sweep per time step loses more accuracy than is gained in computational efficiency by performing fewer sweeps, unless one specifically seeks fast simulation time with less requirements on accuracy. 
However, each additional stage improves accuracy at a rate faster than is added in computational expense, such that the three-stage SSP-LDIRK3(3,3,2) is clearly the most efficient, beating implicit by more than an order of magnitude in wallclock time for fixed accuracy, and as much a two orders of magnitude in terms of accuracy for fixed wall-clock time.
Note that this additional efficiency comes at the expense of storing additional stage vectors in memory. 

\begin{figure}
    \centering
    \begin{subfigure}{.4\textwidth}
        \centering
        \includegraphics[width=\textwidth]{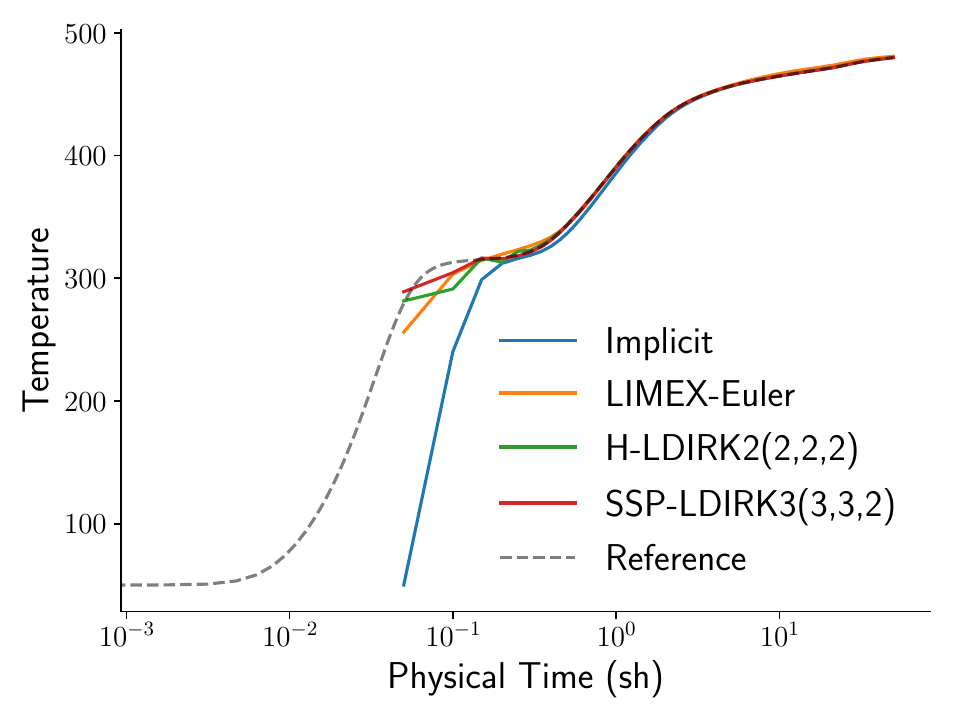}
        \caption{}
        \label{tophat_tracer_a}
    \end{subfigure}
    \begin{subfigure}{.4\textwidth}
        \centering
        \includegraphics[width=\textwidth]{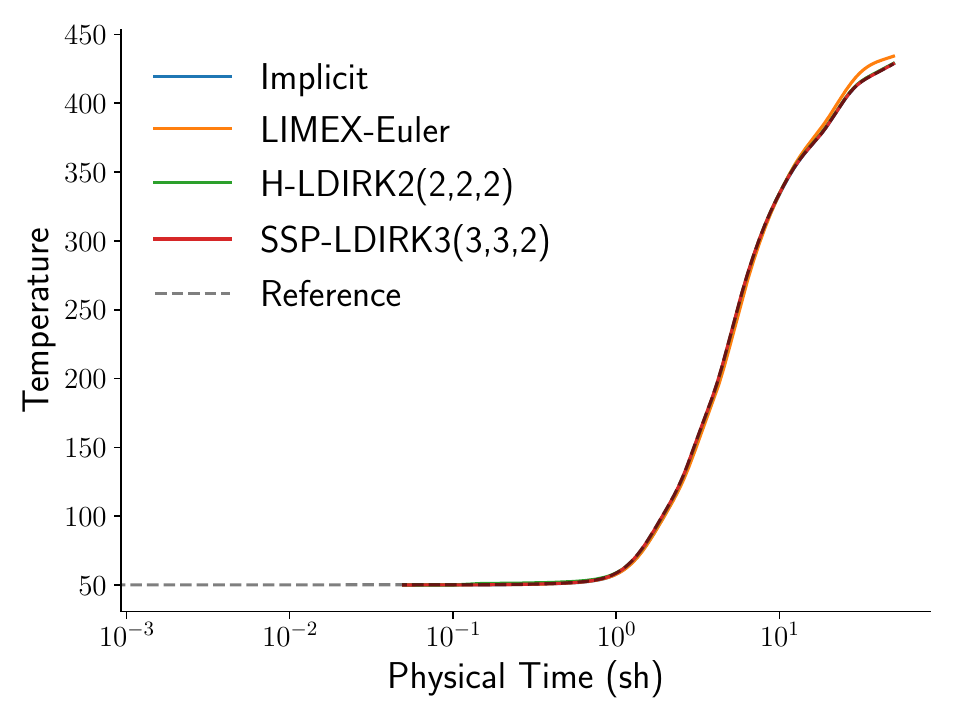}
        \caption{}
    \end{subfigure}
    
    \begin{subfigure}{.4\textwidth}
        \centering
        \includegraphics[width=\textwidth]{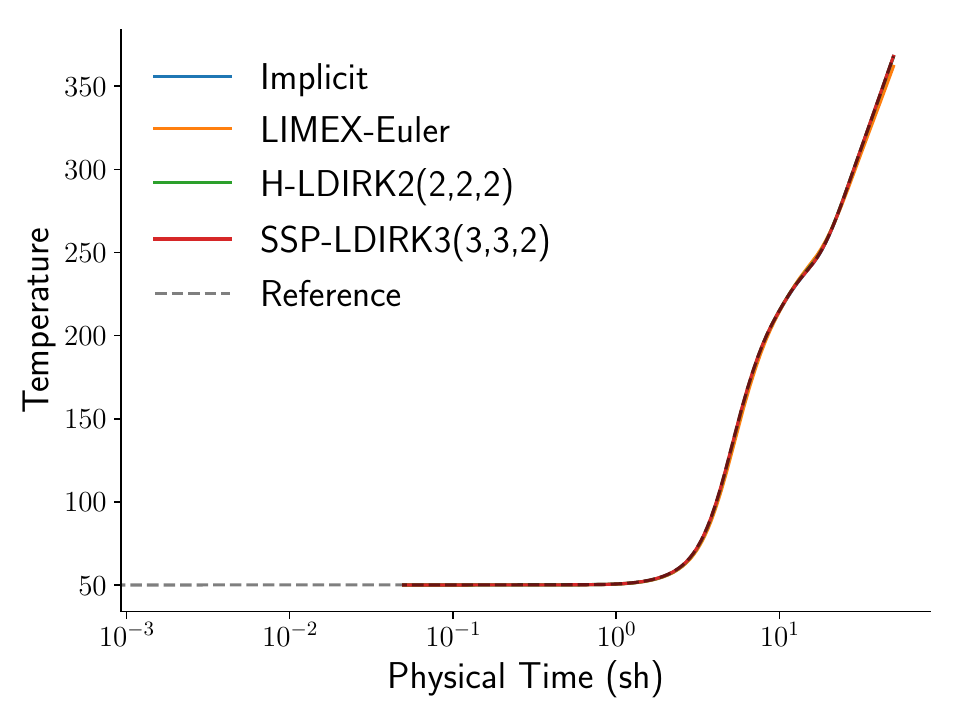}
        \caption{}
    \end{subfigure}
    \begin{subfigure}{.4\textwidth}
        \centering
        \includegraphics[width=\textwidth]{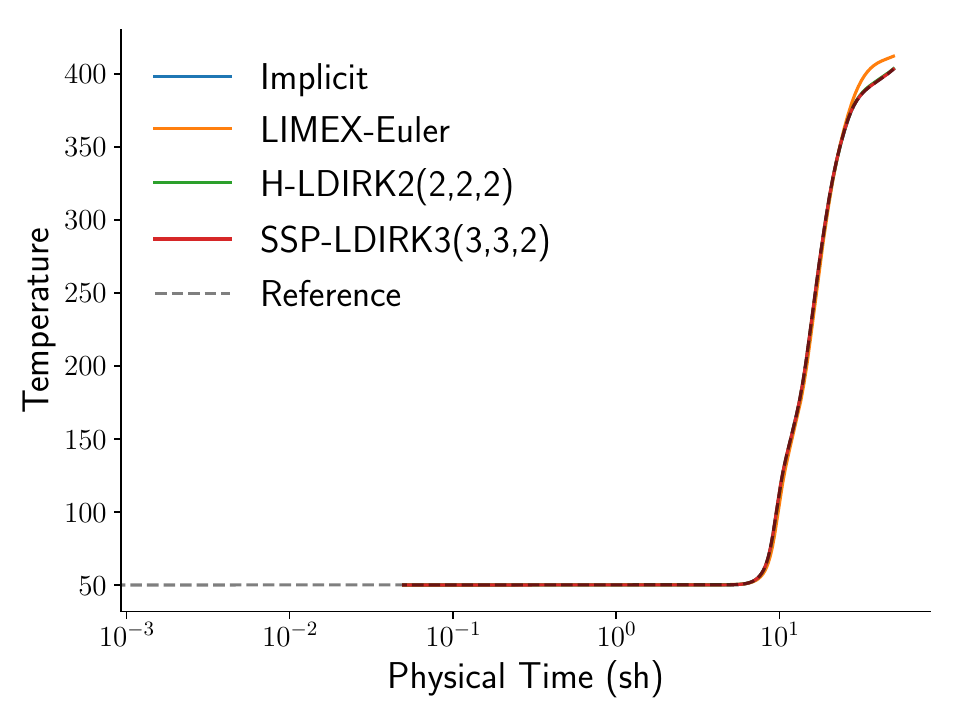}
        \caption{}
    \end{subfigure}
    \caption{Temperature evolution at the four spatial locations depicted in Fig.~\ref{fig:tophat_geometry}. The solution quality generated by the implicit, LIMEX-Euler, H-LDIRK2(2,2,2), and SSP-LDIRK3(3,3,2) HOLO schemes are compared to a reference solution using a time step of size \SI{5e-2}{\sh}. The reference solution is computed with the implicit HOLO scheme using a time step of size \SI{1.5625e-3}{\sh}. }
    \label{fig:tophat_tracer}
\end{figure}
\begin{figure}
\centering
\begin{subfigure}{.4\textwidth}
    \centering
    \includegraphics[width=\textwidth]{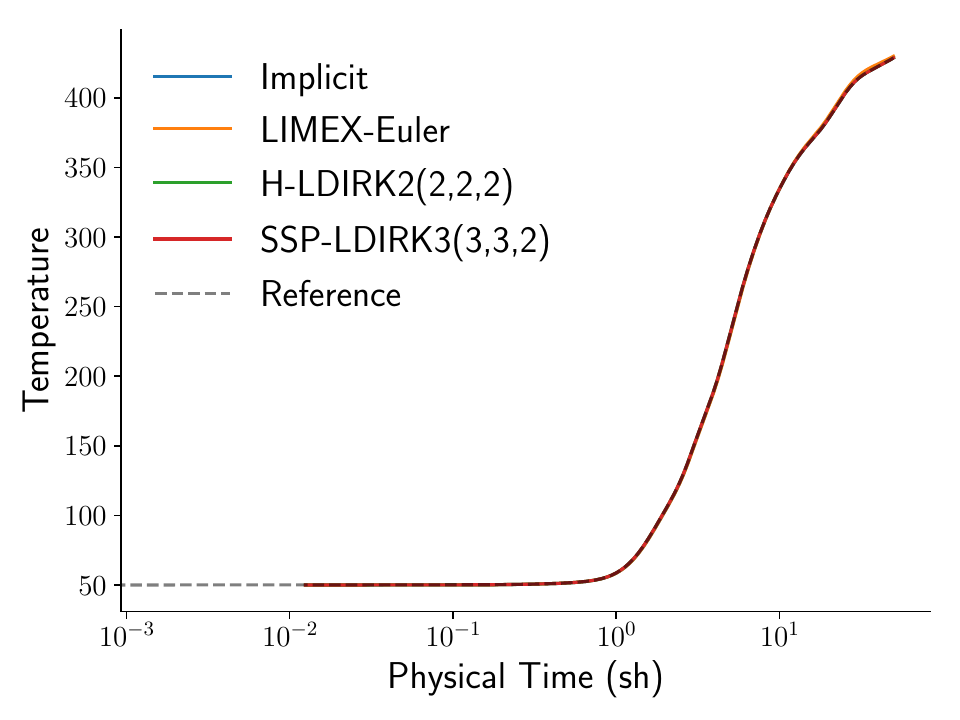}
    \caption{}
\end{subfigure}
\begin{subfigure}{.4\textwidth}
    \centering
    \includegraphics[width=\textwidth]{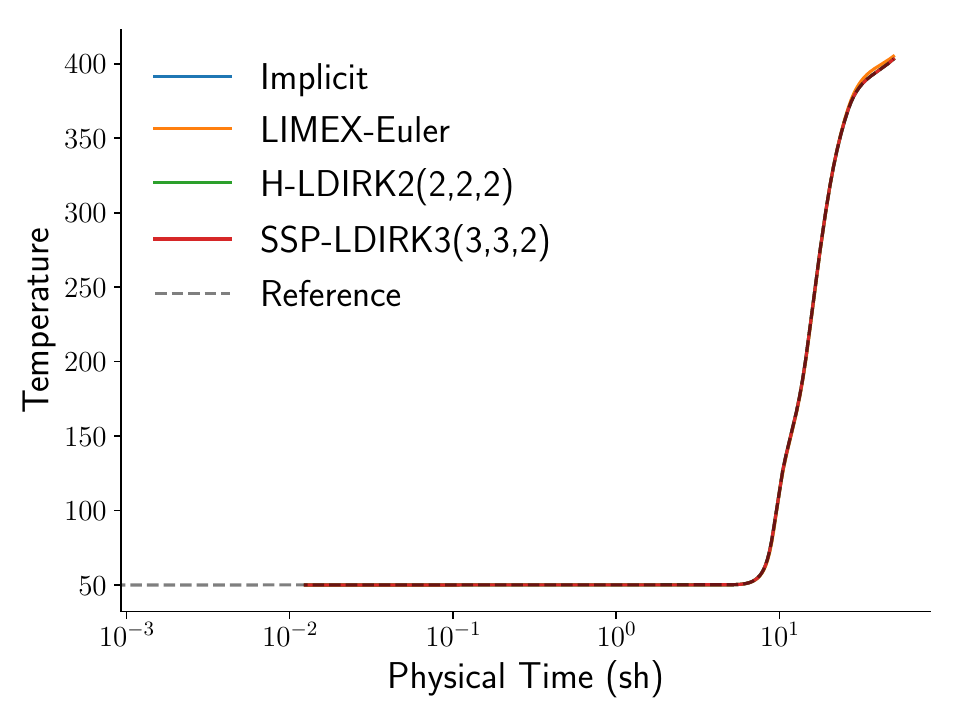}
    \caption{}
\end{subfigure}
\caption{Temperature evolution generated with a time step $4\times$ smaller than in \Cref{fig:tophat_tracer} at locations (b) and (d) where LIMEX-Euler had over predicted the late time behavior of the solution. }
\label{fig:tophat_tracer_refined}
\end{figure}

\Cref{fig:tophat_tracer} shows the time evolution of the temperature at the four spatial locations depicted in \Cref{fig:tophat_geometry} using a time step of size \SI{5e-2}{\sh}. 
We use the implicit HOLO scheme with a time step of size \SI{1.5625e-3}{\sh} as the reference solution. 
All methods track the reference well even at large time steps aside from LIMEX-Euler which shows visually erroneous long time behavior at the mid and endpoints of the pipe. 
Temperature evolution for the pipe mid and endpoints using a time step size $4\times$ smaller is shown in \Cref{fig:tophat_tracer_refined} to show that reducing the time step does improve solution quality for LIMEX-Euler. 

\subsubsection{Comparison to Unaccelerated}
\begin{figure}
    \centering
    \begin{subfigure}{.49\textwidth}
        \centering
        \includegraphics[width=\textwidth]{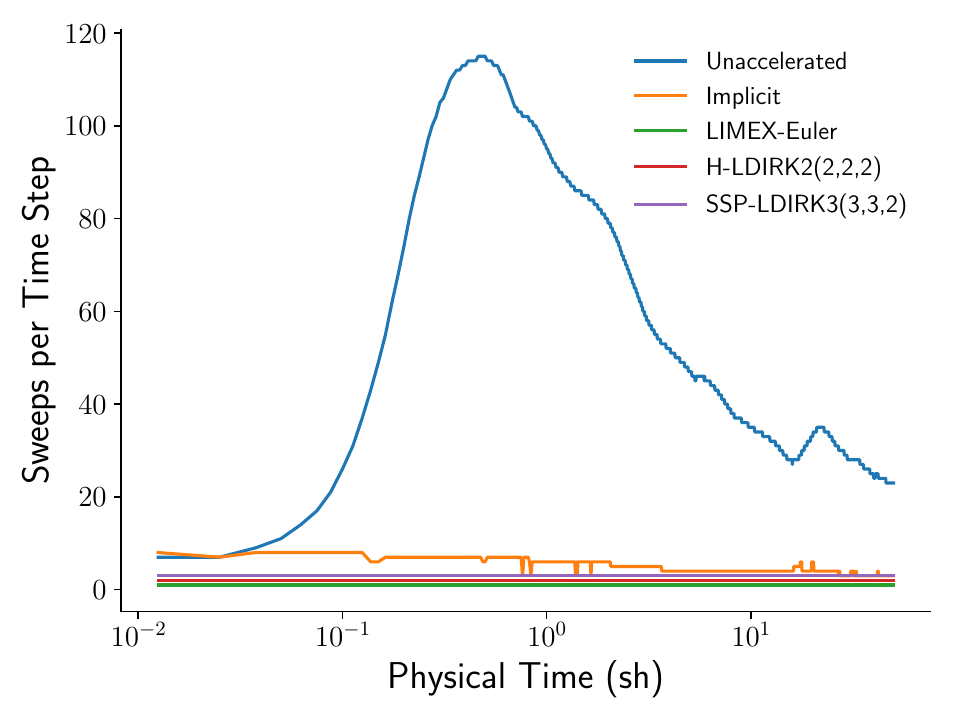}
        \caption{}
    \end{subfigure}
    \begin{subfigure}{.49\textwidth}
        \centering
        \includegraphics[width=\textwidth]{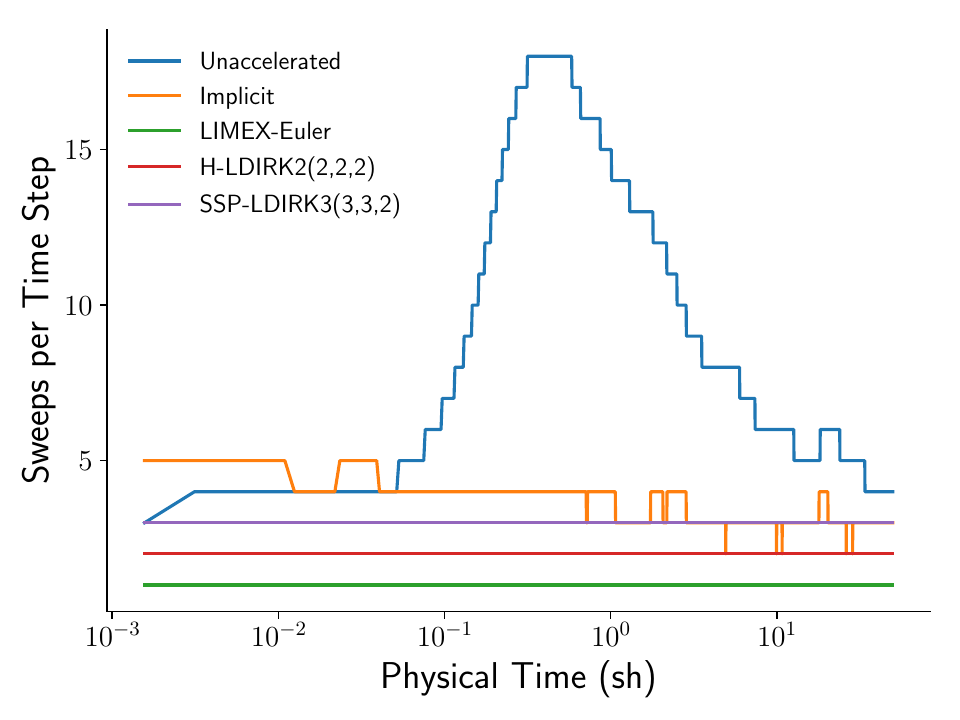}
        \caption{}
    \end{subfigure}
    \caption{Number of sweeps performed at each time step of the crooked pipe problem using a time step of size (a) \SI{1.25e-2}{\sh} and (b) \SI{3.125e-3}{\sh}. An unaccelerated, backward Euler time integration scheme is compared to backward Euler, LIMEX-Euler, H-LDIRK2(2,2,2), and SSP-LDIRK3(3,3,2) HOLO schemes. Note that the LIMEX-Euler, H-LDIRK2(2,2,2), and SSP-LDIRK3(3,3,2) HOLO schemes perform exactly one, two, and three sweeps per time step, respectively, while the implicit schemes require a variable number of sweeps depending on the difficulty of the nonlinear iteration at each time step.  }
    \label{fig:tophat_sweeps}
\end{figure}

We now compare the HOLO schemes to an unaccelerated, backward Euler transport algorithm. 
Here, unaccelerated refers to not using a low-order diffusion system to accelerate the absorption-reemission physics, resulting in an algorithm that requires significantly more iterations to converge when the problem is stiff. 
The time step is set to induce two regimes of performance: one where acceleration is needed to achieve realistic runtimes and one where the time step is small enough such that acceleration is not needed. 
\Cref{fig:tophat_sweeps} shows the number of sweeps performed at each time step when the time step size is \SI{1.25e-2}{\sh} and \SI{3.125e-3}{\sh}. 
The unaccelerated scheme is very sensitive to the solution's dynamics, ranging from performing 7--115 sweeps when the larger time step is used and 3--18 sweeps for the smaller time step. 
By contrast, acceleration via HOLO leads to more predictable cost per time step, with implicit HOLO ranging from 3--8 and 2--5 sweeps per time step on the larger and smaller time step size, respectively.  
By construction, the IMEX schemes perform one sweep per stage so that LIMEX-Euler, H-LDIRK2(2,2,2), and SSP-LDIRK3(3,3,2) perform only one, two, and three sweeps per time step, respectively. 

\begin{figure}
    \centering
    \begin{subfigure}{.49\textwidth}
        \centering
        \includegraphics[width=\textwidth]{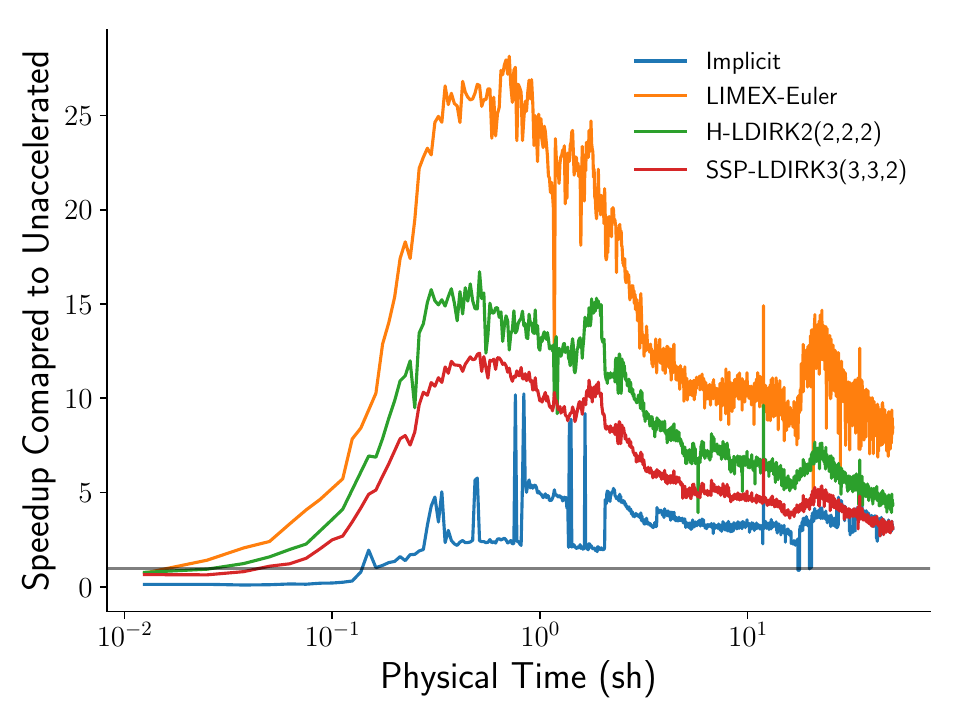}
        \caption{}
    \end{subfigure}
    \begin{subfigure}{.49\textwidth}
        \centering
        \includegraphics[width=\textwidth]{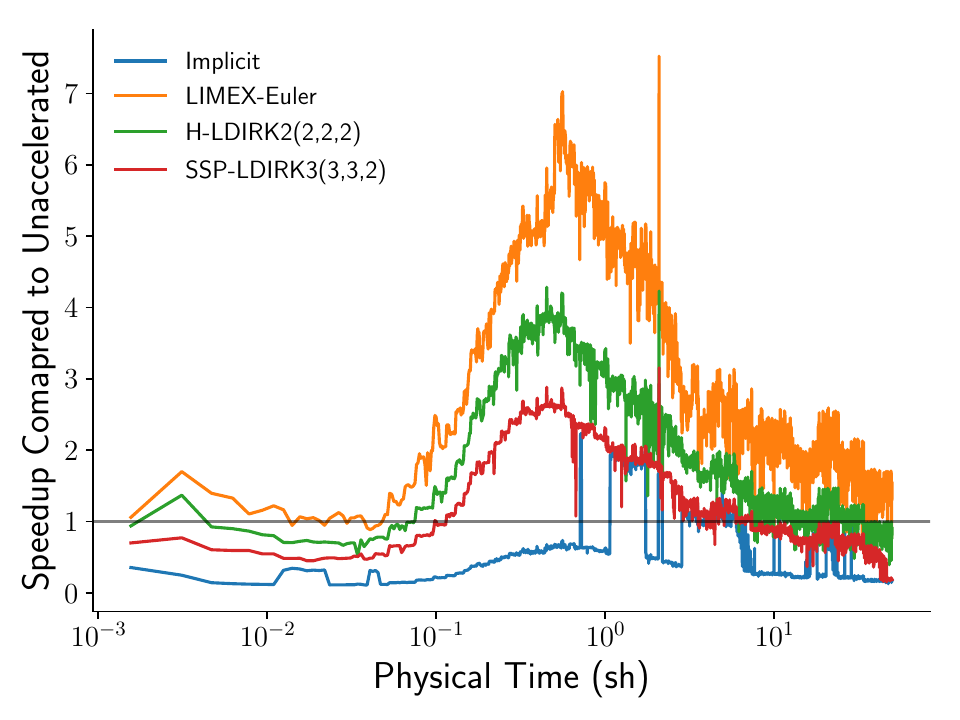}
        \caption{}
    \end{subfigure}
    \caption{Speedup relative to the unaccelerated, backward Euler algorithm in time-to-solution at each time step of the crooked pipe problem using time steps of size (a) \SI{1.25e-2}{\sh} and (b) \SI{3.125e-3}{\sh}. A line at $y=1$ is included to highlight where the schemes cross over from slower than the unaccelerated, backward Euler algorithm to faster. }
    \label{fig:tophat_speedup}
\end{figure}
Speedups compared to the unaccelerated method are shown in \Cref{fig:tophat_speedup} for the large and small time step sizes. 
A horizontal line is included to show the crossover point where the method transitions from slower than unaccelerated to faster. 
For the large time step, the implicit scheme is faster than the unaccelerated scheme after \SI{0.1}{\sh} when the radiation begins to turn the corner. 
At early times, unaccelerated converges rapidly and thus the additional expense of solving the HOLO low-order linear system only serves to slow implicit HOLO down. 
As the problem becomes diffusive, unaccelerated requires more sweeps allowing implicit HOLO to achieve modest speedups of 3--5$\times$. 
On the other hand, the IMEX schemes nearly always perform fewer sweeps than unaccelerated, leading to more consistent speedups over unaccelerated. 
For the larger time step, cost per time step ranges from similar in expense as unaccelerated to 27$\times$, 17$\times$, and 12$\times$ faster for LIMEX-Euler, H-LDIRK2(2,2,2), and SSP-LDIRK3(3,3,2), respectively. 
For the smaller time step, only LIMEX-Euler is always faster than unaccelerated and the maximum speedups are reduced to 7.5$\times$, 4$\times$, and 3$\times$. 

\section{Conclusions}\label{sec:conc}

We have derived a nonlinear moment-based partition of gray TRT to facilitate stable and accurate time integration with a small, fixed number of sweeps per time step. For the thick Marshak problem in 1d, we can match the accuracy of a fully implicit method using only one sweep per time step. For the 2d crooked pipe benchmark, one sweep provides stable time integration, but good accuracy requires schemes with two or three sweeps per time step. Nevertheless, such schemes still provide significant speedups over implicit integration, while also providing significantly improved accuracy, due to the increased order that comes with additional stages as well. The methods provided here are also amenable to tabular opacity data, and can still yield higher than first-order accuracy. Here we only considered a HOLO formulation of gray TRT, but the approach is naturally amenable to other moment formulations as well, such as quasidiffusion/VEF and the second moment method, and future work will extend the methodology to scattering and multifrequency TRT, as well as coupling with hydrodynamics \cite{southworth2023implicit}.  We are also currently deriving a new class of time integration schemes for nonlinear partitions \cite{nprk1,nprk2} such as the one proposed here, and future work will study specialized schemes for efficient integration of TRT and full radiation hydrodynamics.

\section*{Acknowledgements}
This work was supported by the Laboratory Directed Research and Development program of Los Alamos National Laboratory under project number 20220174ER. LA-UR-24-20140. This research used resources provided by the Darwin testbed at Los Alamos National Laboratory (LANL) which is funded by the Computational Systems and Software Environments subprogram of LANL's Advanced Simulation and Computing program (NNSA/DOE).

\bibliographystyle{elsarticle-num} 
\bibliography{refs.bib}
\end{document}